\documentclass[a4paper, 12pt]{article}

\usepackage[swedish, english]{babel}
\usepackage{amsmath}
\usepackage{amssymb}
\usepackage{epsfig}		
\usepackage{float}
\usepackage{subfigure}
\usepackage{moreverb}
\usepackage{epstopdf}

\renewcommand{\appendix}{
\setcounter{section}{0}
\renewcommand{\thesection}{\Roman{section}}
\newpage
\noindent
{\Large{\bf APPENDIX}}}

\newtheorem{theorem}{Theorem}[section]
\newtheorem{lemma}[theorem]{Lemma}

\newtheorem{definition}[theorem]{Definition}

\newenvironment{proof}[1][Proof]{\begin{trivlist}
\item[\hskip \labelsep {\bfseries #1}]}{\end{trivlist}}

\newenvironment{remark}[1][Remark]{\begin{trivlist}
\item[\hskip \labelsep {\bfseries #1}]}{\end{trivlist}}

\newcommand{\qed}{ \nobreak \ifvmode \relax \else
        \ifdim \lastskip<1.5em \hskip-\lastskip
        \hskip1.5em plus0em minus0.5em \fi \nobreak
        \vrule height0.5em width0.5em depth0.25em \fi}

\title{Stable and Accurate Interpolation Operators for High-Order Multi-Block Finite-Difference Methods}
\author{K. Mattsson  \and Mark H. Carpenter}

\date{}

\begin{document}
\selectlanguage{english}
\maketitle
 
\begin{abstract}
Block-to-block interface interpolation operators are constructed for several common high-order finite difference discretizations.
In contrast to conventional interpolation operators, 
these new interpolation operators maintain the strict stability, accuracy
and conservation of the base scheme even when nonconforming grids or dissimilar
operators are used in adjoining blocks.  
The stability properties of the new operators are verified using eigenvalue analysis, 
and the accuracy properties are verified 
using numerical simulations of the Euler equations in two spatial dimensions. 
\end{abstract} 


\section{Introduction}\label{sec:intro}
Wave-propagation problems frequently require farfield boundaries to be positioned many wavelengths away from the disturbance source.  
To efficiently simulate these problems requires numerical techniques capable of accurately propagating disturbances over long
distances. It is well know that high-order finite difference methods (HOFDM) are ideally suited for problems of this type. 
(See the pioneering paper by Kreiss and Oliger \cite{KreissOliger72}). 
Not all high-order spatial operators are applicable, however.  For example, 
schemes that are G-K-S stable \cite{GKS}, while being convergent to the true solution as $\Delta x
\rightarrow 0$, may experience non-physical solution growth in time \cite{CGA}, thereby limiting their efficiency 
for long-time simulations. Thus, it is imperative to use HOFDMs 
that do not allow growth in time; a property termed ``strict stability''\cite{GustafssonKreissOliger}. 
Deriving strictly stable, accurate and conservative HOFDM is a significant challenge and that has received considerable past attention.
(For examples, see references 
\cite{kn:Lele(1992),Strikwerda97,ZinggDeRango99,DeRangoZingg01,Abarbanel97,Bayliss86,Hesthaven98,Sjogreen_95,bertil_pelle2}).

A robust and  well-proven HOFD methodology that ensures the strict stability of 
time-dependent partial differential equations (PDEs) is the summation-by-parts simultaneous approximation term (SBP-SAT) method.  
The SBP-SAT method simply combines 
finite difference operators that satisfy a summation-by-parts (SBP) formula\cite{KreissScherer74}, with physical boundary 
conditions implemented using either the Simultaneous Approximation Term (SAT) method \cite{CGA}, 
or the projection method  \cite{Olsson95(1),Olsson95(2),MattssonNordstrom06}. 
Examples of the SBP-SAT approach can be found in references 
\cite{NordstromCarpenter01,NordstromCarpenter99,Nordstrom_gong_2006,MattssonNordstrom04,MattssonSvard06,MattssonSvard07,NordstromMattsson07,MattssonHam08,SvardNordstrom08a,SvardNordstrom08b,Lehner_05,Diener_07}.

An added benefit of the SBP-SAT method is that it naturally extends to multi-block geometries while retaining the essential
single-block properties: strict stability, accuracy and conservation \cite{CarpenterNordstrom99}.  
Thus, problems involving complex domains or non-smooth geometries are easily amenable to the approach.  
References \cite{MattssonHam08,MattssonSvard06,Svard04_Mattsson} report applications of the 
SBP-SAT HOFD methods to problems involving non-trivial geometries.

Current multi-block SBP-SAT methods suffer from two significant impediments:
1) the collocation points in adjoining blocks must match along block interfaces 
(i.e., ``conforming'' grids must be used), and
2) identical SBP schemes must be used tangential to the block interfaces. 
\footnote{Ironically, these impediments were assumptions needed to prove stability and conservation, 
in the original single domain SBP-SAT approach.} 
These impediments prevent SBP-SAT operators from use on adaptively refined nonconforming multi-block grids and from use
where hybrid-approaches involving two or more discretization techniques is desirable.

Although neither impediment (i.e., conforming grids or identical elements) significantly limits Finite- and Spectral-Elements methods,
\footnote{For example, see the nonconforming spectral multi-domain work of Kopriva and Kolias \cite{kopriva96}.}
to date, mitigating their influence in the context of HOFDM has rarely been successful.  
A noteworthy exception appears in reference \cite{Nordstrom_gong_2006}, where a hybrid method is 
developed (based on the SBP-SAT technique) to merge a HOFDM to an 
unstructured finite-volume method. This was accomplished by alternating the finite volume scheme at the interface boundary. 

The goals of the present study are two-fold.  The first is to develop a systematic methodology for coupling
arbitrary SBP methods between adjoining blocks. 
This task essentially requires identification of the conditions that 
adjoining interpolation operators must satisfy to guarantee interface strict stability, accuracy 
and conservation.

The second task is to demonstrate the new methodology by developing interface interpolation operators that 
couple commonly used SBP operators.  To this end, new interpolation operators are derived to couple 
$2$nd-:$2$nd-, $4$th-:$4$th-, $6$th-:$6$th-, and $8$th-:$8$th-order SBP operators across a 
nonconforming interface containing a $2:1$ grid compression, as well as to couple
$2$nd-:$4$th- and $4$th-:$8$th-order SBP operators across a conforming interface.

In Section \ref{sec:definitions} we introduce some definitions and discuss the SBP property.  
Sufficient conditions for the interface stability of a two-dimensional hyperbolic system are introduced in 
Section \ref{sec:analysis}.  The new approach utilizes 
the SAT technique in combination with newly developed SBP-preserving interpolation operators.
The construction procedure for the SBP-preserving interpolation operators are presented in detail. 
In Section \ref{sec:computations} the accuracy and stability properties of the newly developed multi-block 
interface coupling are tested  by performing an eigenvalue analysis and numerical simulations of an 
analytic Euler vortex. 
In Section \ref{sec:Conclusions} conclusions are drawn. 
Sufficient conditions for interface stability of a 2D parabolic system followed by 
SBP-preserving interpolation operators are presented in the Appendix.

\section{Definitions}\label{sec:definitions}
Two- and three-dimensional schemes are constructed using tensor products of one-dimensional SBP finite-difference operators. 
Thus, we begin with a short description of 1D SBP operators including some relevant definitions.
(For more details see \cite{KreissScherer74,Strand94} and \cite{MattssonNordstrom04}).
\subsection{Summation-By-Parts}\label{sec:1-D}
Assume the existence of real-valued functions $u,v$ such that $u,v \in \mathbf{L^2[0,1]}$. 
Define an inner product to be $(u,v)=\int_0^1u\,v\,dx$, and let the corresponding norm be $\|u\|^2=(u,u)$.  
With these definitions, consider the hyperbolic scalar equation $u_t+u_x=0$ (excluding the boundary
condition). Apply the energy method to $u_t+u_x=0$; i.e., multiplication by $u$ and integration by
parts leads to
\begin{equation}\label{eq:IBP_D1}
\frac{d}{dt}\|u\|^2=-(u,u_x)-(u_x,u)=-u^2|^1_0\;,
\end{equation}
where  $u^2|^1_0 \equiv u^2(x=1)-u^2(x=0)$.  

Next, discretize the domain ($0 \le x \le 1$) using N+1 equidistant grid points,
\[
\begin{array}{lll}
x_i=i\,h,& i=0,1...,N,&  h=\frac{1}{N}\;.
\end{array}
\]
Denote the projection of $v$ onto the discrete grid point $x_i$ as $v_i$, and define
the discrete solution vector to be $v^T=[v_0, v_1, \cdots, v_N]$. 
As with the continuous inner product, 
define the discrete inner product for real-valued vector functions 
$u,v \in \mathbf{R^{N+1}}$ to be $(u,v)_H=u^T\,H\,v$, where $H=H^T\,>\,0$. The 
corresponding discrete norm becomes $\|v\|_H^2=v^T\,H\,v$.  

\begin{definition}\label{definition:SBP_D1a}
A difference operator $D_1=H^{-1}Q$ approximating
$\partial/\partial\,x$ is a first-derivative SBP operator 
if  $H=H^T$, $x^THx\,>\,0$, $x \ne 0$, and $Q+Q^T=B=diag\left(-1,0\ldots,0,1\right)$.
\end{definition}   
A semidiscretization of $u_t+u_x=0$ is $v_t+D_1v=0$. 

A semi-discrete energy for the system can be derived in a manner analogous to that used in the continuous case.
That is, multiplying $v_t+D_1v=0$ from the left by $v^TH$ and adding the transpose leads to,
\begin{equation}\label{eq:SBP_D1}
\frac{d}{dt}\|v\|^2_H=-(v,H^{-1}Qv)_H-(H^{-1}Qv,v)_H=-v^T(Q+Q^T)v=v^2_0-v^2_N\;.
\end{equation}
Equation (\ref{eq:SBP_D1}) is the discrete analog of (\ref{eq:IBP_D1}).

Three definitions are central to the present study.
\begin{definition}\label{definition:narrow}
An explicit p$th$-order accurate finite difference scheme with minimal stencil width of a Cauchy problem, is called a $pth$-order accurate narrow stencil.  
\end{definition}  

\begin{definition}\label{definition:diagonal}
A $pth$-order accurate narrow stencil SBP-operator is called a narrow-diagonal SBP operator if $H$ is diagonal.
\end{definition}  

\begin{definition}\label{definition:accuracy}
Let the  row-vectors ${\bf x^k}_f$ and ${\bf x^k}_c$ be the projections of the polynomials $x^k$ onto equidistant 1-D grids  corresponding to a fine and coarse grid respectively.  We say that $I_{F2C}$ and $I_{C2F}$ are p$th$-order accurate interpolation operators  if $e^k_c \equiv I_{F2C}{\bf x^k}_f-{\bf x^k}_c$ and $e^k_f \equiv I_{C2F}{\bf x^k}_c-{\bf x^k}_f$ vanish for $k=0..p-1$ in the interior and for $k=0 \ldots (p-1)/2$ at the boundaries.
\end{definition}

\begin{remark}
The boundary closure for a $pth$-order accurate narrow-diagonal SBP operator is of order $p/2$  (see \cite{MattssonNordstrom04}).  
The convergence rate for narrow stencil approximations of fully hyperbolic problems (e.g., the Eulers equations)
drops to $(p/2+1)th$-order.  (See refs. \cite{Gustafsson81,SvardNordstrom06} for more information on the accuracy of 
finite difference approximations).
\end{remark}

\begin{remark}
We state here without proof that 
the accuracy of the $p$th-order narrow-diagonal SBP operator is preserved when
using the $p$th-order accurate interpolation operators.  Numerical experiments
presented later herein, are consistent with this conjecture.
\end{remark}

\subsection{Two-Dimensional Domains}\label{sec:2-D}
We begin by introducing the Kronecker product
\[
C \otimes D=\begin{bmatrix}
c_{0,0}\,D &\cdots  &c_{0,q-1}\,D\\
\vdots &  &  \vdots\\
c_{p-1,0}\,D &\cdots  &c_{p-1,q-1}\,D\\
\end{bmatrix}\;,
\]
where $C$ is a $p \times q$ matrix and $D$ is an $m \times n$ matrix. Two useful
rules for the Kronecker product are $(A \otimes B)(C \otimes D)=(AC)
\otimes (BD)$ and $(A \otimes B)^T=A^T \otimes B^T$.

\begin{figure}
\centering
\epsfig{file=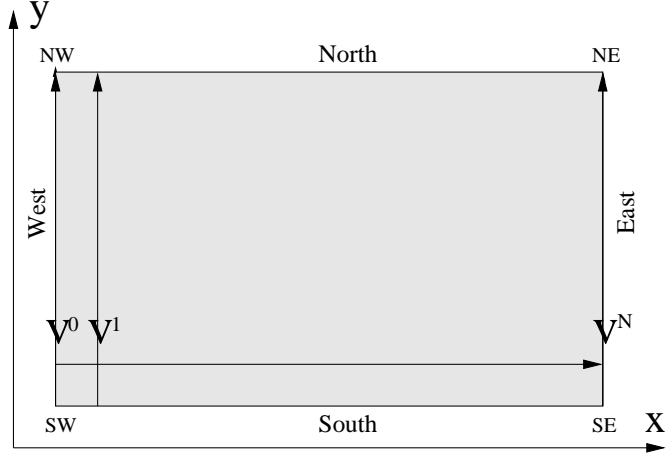,height=6cm,angle=0}
\caption{Domain 2-D}
\label{fig:domain_2d}
\end{figure}

Next, consider the 2-D rectangular domain  $\Omega^{l,\,r}_{d,\,u}$ 
bounded by $l\le x \le r,\;d\le y \le u$ with an $(N+1)\times (M+1) $-point equidistant grid defined as,
\[
\begin{array}{lll}
x_i=ih_x,&i=0,1...,N,& h_x=\frac{r-l}{N},\\
y_j=jh_y,& j=0,1...,M,& h_y=\frac{u-d}{M}.\\
\end{array}\;
\]
The numerical approximation at grid point $(x_i,\,y_j)$ is a $1\times k$-vector denoted $v_{i,j}$.
The tensor product derivations are more transparent if we redefine the component vector $v_{i,j}$ as a ``vector of vectors''.  Specifically,
define a discrete solution vector $v^T=[v^0,v^1,\ldots,v^N]$, where $v^p=[v_{p,0}, v_{p,1}, \ldots, v_{p,M}]$ is the solution vector at $x_p$ along the y-direction, illustrated in Figure \ref{fig:domain_2d}.
To distinguish whether a difference operator $P$  is operating in the $x-$ or the $y$-direction 
we use the notations $P^x$ and $P^y$. 
The following two-dimensional operators are frequently used,
\begin{equation}\label{eq:Def_matrix2}
\begin{array}{llll}
 {\bf D^x}&= (D^x \otimes I^y),& {\bf D^y}&= (I^x \otimes D^y)\\
 {\bf H^x}&= (H^x \otimes I^y),& {\bf H^y}&= (I^x \otimes H^y)\\
 {\bf e^x_0}&= (e_0 \otimes I^y),& {\bf e^y_0}&= ( I^x \otimes e_0 )\\
  {\bf e^x_N}&= (e_N \otimes I^y),& {\bf e^y_M}&= ( I^x \otimes e_M )\\
  \end{array}\;,
\end{equation}
where $D^{x,\,y}$, and $H^{x,\,y}$ are the one-dimensional operators in the  x- and y-direction, respectively, 
$I^{x,\,y}$ are identity matrices of appropriate sizes, and $e_0$, and $e_N$ are one-dimensional 
``boundary'' vectors defined by
\begin{equation}\label{eq:Def_matrix}
e_0 = [1,0,...,0]^T,\;e_N = [0,...,0,1]^T\;\:.
\end{equation}

\section{Analysis}\label{sec:analysis}
Our main focus is on multi-block interface coupling solving PDEs on nonconforming grids, combining different (here referring to the formal order of accuracy) finite difference SBP schemes. 

Consider the following 2-D hyperbolic problem (the extension to parabolic problems is shown in the Appendix):
\begin{equation}\label{eq:hyper_pde}
\begin{array}{ll}
\medskip {\bf u}_t+A{\bf u}_x+B{\bf u}_y=0,&{\bf x} \in \Omega^{l,\,i}_{d,\,u},\\
\medskip {\bf v}_t+A{\bf v}_x+B{\bf v}_y=0,&{\bf x} \in \Omega^{i,\,r}_{d,\,u}\\
\end{array}\;,
\end{equation}
where $A$ and $B$ are symmetric $k\times k$ matrices ($k=4$ for the compressible Euler equations). (The extension of the present study to 3-D is straightforward.) A corresponding semi-discrete approximation of (\ref{eq:hyper_pde}) can be written
\begin{equation}\label{eq:hyper_SAT}
\begin{array}{l}
\medskip u_t+A \otimes {\bf D_L^x}u+B \otimes {\bf D_L^y}u=SAT^I_{IL}+SAT^I_L\\
\medskip v_t+A \otimes {\bf D_R^x}v+B \otimes {\bf D_R^y}v=SAT^I_{IR}+SAT^I_R\;,\\
\end{array}
\end{equation}
where $SAT^I_{L,\,R}$ denote the penalty terms at the outer boundaries. (These are tuned to obtain stability, see \cite{NordstromCarpenter01,MattssonSvard07,Svard04_Mattsson} for details.) 
The penalty terms  handling the multi-block coupling can be written as 
\begin{equation}\label{eq:hyper_stability}
\begin{array}{l}
\medskip SAT^I_{IL}=\Sigma_L \otimes ({\bf H^{x}_L)^{-1}}{\bf e^x_N}\left\{u_N-I_k\otimes I_{F2C}v_0\right\}\\
\medskip SAT^I_{IR}=\Sigma_R \otimes ({\bf H^{x}_R})^{-1}{\bf e^x_0}\left\{v_0-I_k\otimes I_{C2F}u_N\right\}\\
\end{array}\;.
\end{equation}  
(However, the penalty terms (\ref{eq:hyper_stability}) is not the most general representation, as will be shown later in this section.)
The continuity condition $u-v=0$ along the interface (at $x_i$) is approximated by $u_N-I_k\otimes I_{F2C}v_0$ and $v_0-I_k\otimes I_{C2F}u_N$ in the above penalty-terms.
$\Sigma_R$ and $\Sigma_{L}$ are unknown $k \times k$ matrices to be determined below by stability.
We use the subscripts $L,\,R$ to distinguish the operators in the left and right domain (which can be completely unrelated, i.e., we can use different numerical schemes in the two domains), with the only assumption that they fulfill the SBP property. $I_{F2C}$ and $I_{C2F}$ are 1-D interpolation operators constructed such that $I_k\otimes I_{F2C}v_0$ and   $I_k\otimes I_{C2F}u_N$ have the same dimension as $u_N$ and $v_0$ respectively. $u_N$ and $v_0$ are $km\times 1$ vectors with components aligned with the block interface. For future reference we introduce the following norms:  
\begin{equation}\label{eq:norms}
 M_L=(I_k \otimes H_L^x\otimes H_L^y)\;,\;M_R=(I_k \otimes H_R^x\otimes H_R^y)\;,
 \end{equation}
and the following terms:
\begin{equation}\label{eq:hyper_terms}
IT^I=
\begin{array}{ll}
\medskip -u^T_N(A\otimes H^y_L)u_N&+v^T_0(A\otimes H^y_R)v_0\\
\medskip +2u^T_N(\Sigma_{L}\otimes H^y_L)u_N&+2v^T_0(\Sigma_R\otimes H^y_R)v_0\\
\medskip -u^T_N(\Sigma_{L}\otimes H^y_LI_{F2C})v_0&-v^T_0(\Sigma_R\otimes H^y_RI_{C2F})u_N\\
\medskip -u^T_N(\Sigma_R\otimes I_{C2F}^TH^y_R)v_0 &-v^T_0(\Sigma_{L}\otimes I_{F2C}^TH^y_L)u_N\\
\end{array}\;,
\end{equation}
\begin{equation}\label{eq:matrix_SAT}
X=
\begin{bmatrix}
 H^y_L               & - H^y_L I_{F2C} \\
-H^y_R I_{C2F}                  & H^y_R\\
\end{bmatrix}\;,
w^T=
\begin{bmatrix}
u_{N}\\
v_0\\
\end{bmatrix}\;.
\end{equation}

We introduce the following important relations,
\begin{eqnarray}
\medskip H^y_RI_{C2F}=I_{F2C}^TH^y_L \label{inter1},\\
\medskip H^y_L(I^y_L-I_{F2C}I_{C2F})\ge0 \;,\; H^y_R(I^y_R-I_{C2F}I_{F2C})\ge0  \label{inter2},
\end{eqnarray}
where $I^y_{L,\,R}$ are  identity matrices of appropriate sizes. 
The following definition is central to the present study, 
\begin{definition}\label{definition:interpolation}
A p$th$-order accurate interpolation operator that fulfills (\ref{inter1}) is called a p$th$-order accurate SBP-preserving interpolation operator.  
\end{definition}

The following lemma is central to the present study:
\begin{lemma}\label{lemma:lemma1}
The scheme (\ref{eq:hyper_SAT}) with the penalty terms given by (\ref{eq:hyper_stability}) is stable if  
$2 \Sigma_{L}=A \;,\; 2 \Sigma_R=-A$ and (\ref{inter1}) hold,
assuming that the boundary terms $SAT^I_{L,\,R}$ are correctly implemented.
\end{lemma}

\begin{proof}
Apply the energy method by multiplying the first and second equation in (\ref{eq:hyper_SAT}) by $u^TM_L$ and $v^TM_R$  (defined by (\ref{eq:norms}))  respectively, yielding,
\[
 \frac{d}{dt}(\|u\|_{M_L}+\|v\|_{M_R})=BT^I+IT^I\,,
\]
where $BT^I$ corresponds to the outer boundary terms (which we assume are correctly implemented, i.e., bounded) and $IT^I$ corresponds to the multi-block interface coupling given by (\ref{eq:hyper_terms}).

A  stable coupling is obtained if $IT^I$ is non-positive. A non-dissipative and stable coupling is achieved by choosing
\[
2 \Sigma_{L}=A \;,\; 2 \Sigma_R=-A \;,\;H^y_RI_{C2F}=I_{F2C}^TH^y_L\;,
\]
which makes $IT^I=0$.
 $\square$
\end{proof}
\begin{remark}
Note the critical relationship that equation (\ref{inter1}) 
imposes on the two interface interpolation operators.
Indeed, the two conditions $2 \Sigma_{L}=A \;,\; 2 \Sigma_R=-A$ and equation (\ref{inter1}) provide 
sufficient conditions for interface stability.
\end{remark}
\begin{remark}
Another important interface property is conservation, especially for problems with discontinous solutions.
Indeed, the Lax-Wendroff theorem relates discrete conservation to convergence to the weak
solution of the differential equation.  It is shown in \cite{CarpenterNordstrom99,NordstromCarpenter99} for
conforming grids that conservation provides a necessary condition for interface stability. 
Likewise, the conditions for conservation are fulfilled by imposing the necessary conditions 
for stability in Lemma \ref{lemma:lemma1} (or Lemma \ref{lemma:lemma2}). 
\end{remark}
For non-linear problems (such as the compressible Euler equations) characteristic boundary conditions are often imposed (since they introduce damping \cite{MattssonSvard06} at the boundaries)  which can also be utilized at block interfaces (see for example \cite{MattssonSvard06,Svard04_Mattsson}). 
 Since $A$ is a symmetric matrix, it can be diagonalized by $R^TAR=\Lambda$, where $R$
is a matrix consisting of the eigenvectors of $A$. 
We introduce the notation $\bar{A} \equiv R|\Lambda|R^T$, and $2A^\pm \equiv A \pm \bar{A}$. 
(The matrix $B$ is also symmetric and can be rotated to diagonal form with a similar transformation.) 
To couple characteristic variables using the SAT technique \cite{MattssonSvard06,Svard04_Mattsson} at the block-interface we choose $\Sigma_{L}=A^- \;,\;  \Sigma_R=-A^+$. 

The following lemma states the stability conditions for a characteristic multi-block coupling: 
\begin{lemma}\label{lemma:lemma2}
The scheme (\ref{eq:hyper_SAT}) with the penalty terms given by (\ref{eq:hyper_stability}) is stable if  $ \Sigma_{L}=A^- \;,\;  \Sigma_R=-A^+ $, (\ref{inter1}) and  (\ref{inter2}) hold,
assuming that the boundary terms $SAT^I_{L,\,R}$ are correctly implemented. 
\end{lemma} 

\begin{proof}
Apply the energy method by multiplying the first and second equation in (\ref{eq:hyper_SAT}) by $u^TM_L$ and $v^TM_R$  (defined by (\ref{eq:norms}))  respectively, yielding,
\[
 \frac{d}{dt}(\|u\|_{M_L}+\|v\|_{M_R})=BT^I+IT^I\,,
\]
where $BT^I$ corresponds to the outer boundary terms (which we assume are correctly implemented, i.e., bounded) and  $IT^I$ corresponds to the multi-block interface coupling given by (\ref{eq:hyper_terms}).

A  stable coupling is obtained if $IT^I$ is non-positive. If $\Sigma_{L}=A^- \;,\;  \Sigma_R=-A^+ $and $H^y_RI_{C2F}=I_{F2C}^TH^y_L$ hold, $IT^I=-\bar{A} \otimes w^TXw$ where $w^T$ and $X$ are given by (\ref{eq:matrix_SAT}).
(By assumption $H^y_L I_{F2C}=(H^y_R I_{C2F})^T$ which makes $X$ symmetric.)
Hence, stability follows if $X$ is positive semi-definite, since $\bar{A}$ is a symmetric and positive semi-definite matrix. Sylvester's theorem states that $X\ge0$ if the following conditions hold:
\[
\begin{array}{ll}
\medskip H^y_L > 0, & H^y_L - H^y_L I_{F2C} (H^y_R)^{-1} H^y_R I_{C2F}\ge0\\
\medskip H^y_R > 0, & H^y_R-  H^y_R I_{C2F}(H^y_L)^{-1}H^y_L I_{F2C}\ge0\\
\end{array}\;.
\]
By definition $H^y_L > 0$, $H^y_L > 0$.  Hence stability follows if (\ref{inter2}) holds.
$\square$
\end{proof}

The extra conditions in (\ref{inter2}), besides the necessary conditions in (\ref{inter1}) introduce a significant restriction in the construction  of the interpolation operators (as will be shown in the coming section). The energy estimate for the characteristic coupling  (as described in Lemma \ref{lemma:lemma2}) gave us an idea how to remove the extra conditions (\ref{inter2}) and yet introduce a damping mechanism by the penalty coupling. Consider the following energy estimate
\begin{equation}\label{eq:norms2}
 \frac{d}{dt}(\|u\|_{M_L}+\|v\|_{M_R})=BT^I-\Omega \otimes w^TX^{\alpha} w\,,
\end{equation}
where $\Omega$ is symmetric and positive definite and $\alpha$ is a positive integer ($w^T$ and $X$ are given by (\ref{eq:matrix_SAT})). For  $\Omega= \bar{A}$ and $\alpha =1$ we recover the energy estimate from the characteristic coupling,  presented  in Lemma \ref{lemma:lemma2}. By choosing $\alpha$ an even integer, stability follows if $X$ is symmetric (which is true if condition (\ref{inter1}) holds),  regardless of the extra conditions in (\ref{inter2}). 
We introduce the following penalties:
\begin{equation}\label{eq:hyper_stability2}
\begin{split}
\medskip SAT^I_{IL}=&+\frac{A}{2} \otimes ({\bf H^{x}_L)^{-1}}{\bf e^x_N}\left\{u_N-I_k\otimes I_{F2C}v_0\right\}\\
\medskip &-\Omega \otimes ( {\bf H^{x}_L)^{-1}  H^{y}_L }{\bf e^x_N}\left\{u_N-I_k\otimes I_{F2C}v_0\right\}\\
\medskip&-\Omega \otimes ( {\bf H^{x}_L)^{-1} I_{F2C}  H^{y}_R }{\bf e^x_N}\left\{I_{C2F}u_N- v_0\right\}\\
\medskip SAT^I_{IR}=&-\frac{A}{2} \otimes ({\bf H^{x}_R})^{-1}{\bf e^x_0}\left\{v_0-I_k\otimes I_{C2F}u_N\right\}\\
\medskip &-\Omega \otimes ( {\bf H^{x}_R)^{-1}  H^{y}_R }{\bf e^x_0}\left\{v_0-I_k\otimes I_{C2F}u_N\right\}\\
\medskip&-\Omega \otimes ( {\bf H^{x}_R)^{-1} I_{C2F}  H^{y}_L }{\bf e^x_0}\left\{I_{F2C}v_0- u_N\right\}\\
\end{split}\;.
\end{equation}  
The first two penalty terms in (\ref{eq:hyper_stability2}) (,i.e., the first in $SAT^I_{IL}$ and the first in $SAT^I_{IR}$ ) are identical to the non-dissipative penalty coupling defined in Lemma \ref{lemma:lemma1}. The following lemma resolves the issue with the extra conditions in (\ref{inter2}),
\begin{lemma}\label{lemma:lemma3}
The scheme (\ref{eq:hyper_SAT}) with the penalty terms given by (\ref{eq:hyper_stability2})
is stable if $\Omega$ is symmetric and positive definite and  (\ref{inter1}) holds,
assuming that the boundary terms $SAT^I_{L,\,R}$ are correctly implemented. 
\end{lemma} 

\begin{proof}
Apply the energy method by multiplying the first and second equation in (\ref{eq:hyper_SAT}) by $u^TM_L$ and $v^TM_R$  (defined by (\ref{eq:norms}))  respectively, leading to  (\ref{eq:norms2}), with $\alpha \equiv 2$. $BT^I$ corresponds to the outer boundary terms (which we assume are correctly implemented, i.e., bounded). Stability follows if $\Omega$ is symmetric and positive definite.
 $\square$
\end{proof}
Lemma \ref{lemma:lemma3} is validated numerically  in Appendix \ref{sec:energy2}, where we also show that condition (\ref{inter1}) is necessary to guarantee a stable characteristic coupling.

\begin{remark}
For non-linear phenomena the addition of artificial dissipation (AD) is often required for stability. The addition of SBP preserving AD \cite{MattssonSvard04} does not alter the stability conditions in Lemma \ref{lemma:lemma1}, Lemma \ref{lemma:lemma2} or Lemma \ref{lemma:lemma3}. 
\end{remark}

\subsection{Construction of Interpolation Operators}\label{sec:construction}
This section briefly describes the construction of the SBP-preserving  interpolation operators $I_{F2C}$ and $I_{C2F}$. In the present study we will consider a two to one ratio between the coarse and the fine grids along the common interface. (The analysis is not restricted to a two to one ratio, it merely simplifies the construction.)
The structure (here showing the upper part) of the interpolation operator $I_{F2C}$ is given by:
\[
\left(
\begin{array}{lllllllllllll}
a_{1,1}&&&\cdots&&\cdots&&\cdots&&&a_{1,r}\\
\vdots&&&&&&&&&&\vdots\\
a_{q,1}&&&\cdots&&\cdots&&\cdots&&&a_{q,r}\\
a_{2s}&a_{2s-1}&a_{2s-2}&\cdots&a_0&\cdots&a_{2s-2}&a_{2s-1}&a_{2s}&&\\
 & &a_{2s}&a_{2s-1}&a_{2s-2}&\cdots&a_0&\cdots&a_{2s-2}&a_{2s-1}&a_{2s}\\\\
&&&\ddots&&\ddots&&\ddots& &\ddots& \ddots\\
\end{array}
\right)
\]
(The lower part of $I_{F2C}$ is obtained by a permutation of both rows and columns.)
The other interpolation operator is given by$I_{C2F}=(H_L^y)^{-1}(H_R^yI_{F2C})^T$. 
For example, in the fourth-order case (see the Appendix) we have used $q=3, r=11$ and $s=2$ , which leads to the following form:
\[
I_{C2F}=
\left(
\begin{array}{lllllllllll}
b_{1,1}&&&\cdots&&&&b_{1,8}&&&\\
\vdots&&&&&&&\vdots&&&\\
b_{11,1}&&&\cdots&&&&b_{11,8}&&&\\
 & & & &2a_3&2a_1&2a_1&2a_3&&&\\
 & & & &2a_4&2a_2&2a_0&2a_2&2a_4&&\\
 & & & & & &2a_3&2a_1&2a_1&2a_3&\\
 & & & & & &2a_4&2a_2&2a_0&2a_2&2a_4\\
&&&&&&\ddots&&\ddots&&\ddots\\
\end{array}
\right)
\]
(The lower part of $I_{C2F}$ is obtained by a permutation of both rows and columns.)
Notice that the unknowns in $I_{C2F}$ are completely defined by the unknowns in $I_{F2C}$ and the known coefficients in $H_L^y$ and $H_R^y$.
\begin{remark}
A p$th$-order accurate interpolation operator will preserve the order of the finite difference scheme (\ref{eq:hyper_SAT}) utilizing narrow-diagonal SBP operators, i.e., resulting in (p/2+1)$th$ order accuracy. 
SBP discretizations not restricted  (see for example\cite{MattssonNordstrom04,Strand94})  to  a $(p/2)th$-order accurate  boundary stencil, require a more accurate boundary closure of the interpolation operators  to maintain the convergence rate.  We omit a detailed study (see  \cite{SvardNordstrom06} for more information on the accuracy of finite difference approximations) of the accuracy requirement. 
\end{remark}

The following procedure describes the construction of the $pth$-order accurate interpolation operators, given the y-norms in the left ($H^y_L$) and right ($H^y_R$) domains:
\begin{itemize}
\item Build  $I_{F2C}$ with the form given above.
\item Build $I_{C2F}=(H^y_R)^{-1}I_{F2C}^TH^y_L$, given $H^y_L$ and $H^y_R$.
\item Solve for accuracy (see Definition \ref{definition:accuracy}).
\end{itemize}
(In step 2 above we have assumed that the left domain corresponds to the coarse grid, if not simply replace the positions between $H^y_L$ and $H^y_R$ if we have the opposite scenario.)  
If there are  free parameters left after imposing the accuracy conditions they are tuned to obtain an optimal $L_2$ error. For the  p$th$-order accurate case this means that $e_k^2 \equiv (e^k_{c})^T \cdot e^k_{c}+(e^k_{f})^T \cdot e^k_{f}$ (see Definition \ref{definition:accuracy}) is optimized for $k=(p-1)/2+1 \ldots$ (resulting in linear equations to be solved) until there are no free parameters left.

With these assumptions there exist by construction SBP-preserving  interpolation operators for the second, fourth, sixth and eighth  order cases (based on narrow-diagonal SBP schemes), by using the symbolic mathematics software Maple, (see the Appendix). We have also constructed interpolation operators coupling fourth to second and eighth to fourth-order accurate SBP discretizations. 
\begin{remark}
This is the first time (to our knowledge) that stability for multi-block coupling of finite difference schemes with different order of accuracy in more than 1-D have been shown. 
\end{remark}
For the second-order and fourth-order cases the constructed interpolation operators fulfill the extra conditions in (\ref{inter2}), besides the necessary conditions in (\ref{inter1}).  
It is an open question if the conditions in (\ref{inter2}) can be met for the sixth and eighth-order cases.

\section{Computations}\label{sec:computations}
We perform some numerical tests to verify the accuracy and stability properties of the SBP-preserving interpolation operators. 

\subsection{Numerical Validation of Lemma \ref{lemma:lemma1}}\label{sec:energy}
To test the stability properties numerically we can either perform : 1) a long-time integration, or 2) an eigenvalue analysis. We chose the latter in the present study.
Consider (\ref{eq:hyper_SAT}) with the following setup:
\begin{equation}\label{eq:PDE_S2}
\begin{array}{lll}
A=B=\begin{bmatrix}
1&0\\
0&-1\\
\end{bmatrix},
&
u=\begin{bmatrix}
u^{(0)}\\
u^{(1)}\\
\end{bmatrix}
& 
v=\begin{bmatrix}
v^{(0)}\\
v^{(1)}\\
\end{bmatrix}\\
u^{(0)}=u^{(1)},&v^{(0)}=v^{(1)}:& \text{at the outer boundaries}\\
u^{(0)}=v^{(0)},& u^{(1)}=v^{(1)}:& \text{at the block-interface} \\
\end{array}\;.
\end{equation}
Despite it's simplicity, this model problem  is very challenging to solve numerically (without introduction of artificial damping) since it is energy-conserving, i.e., the eigenvalues are purely imaginary (see \cite{Mattsson03} for details)).  This makes it an ideal benchmark problem testing for stability of a certain numerical scheme, including boundaries. 

\begin{figure}
  \subfigure{\includegraphics[height=6.5cm,angle=-90]{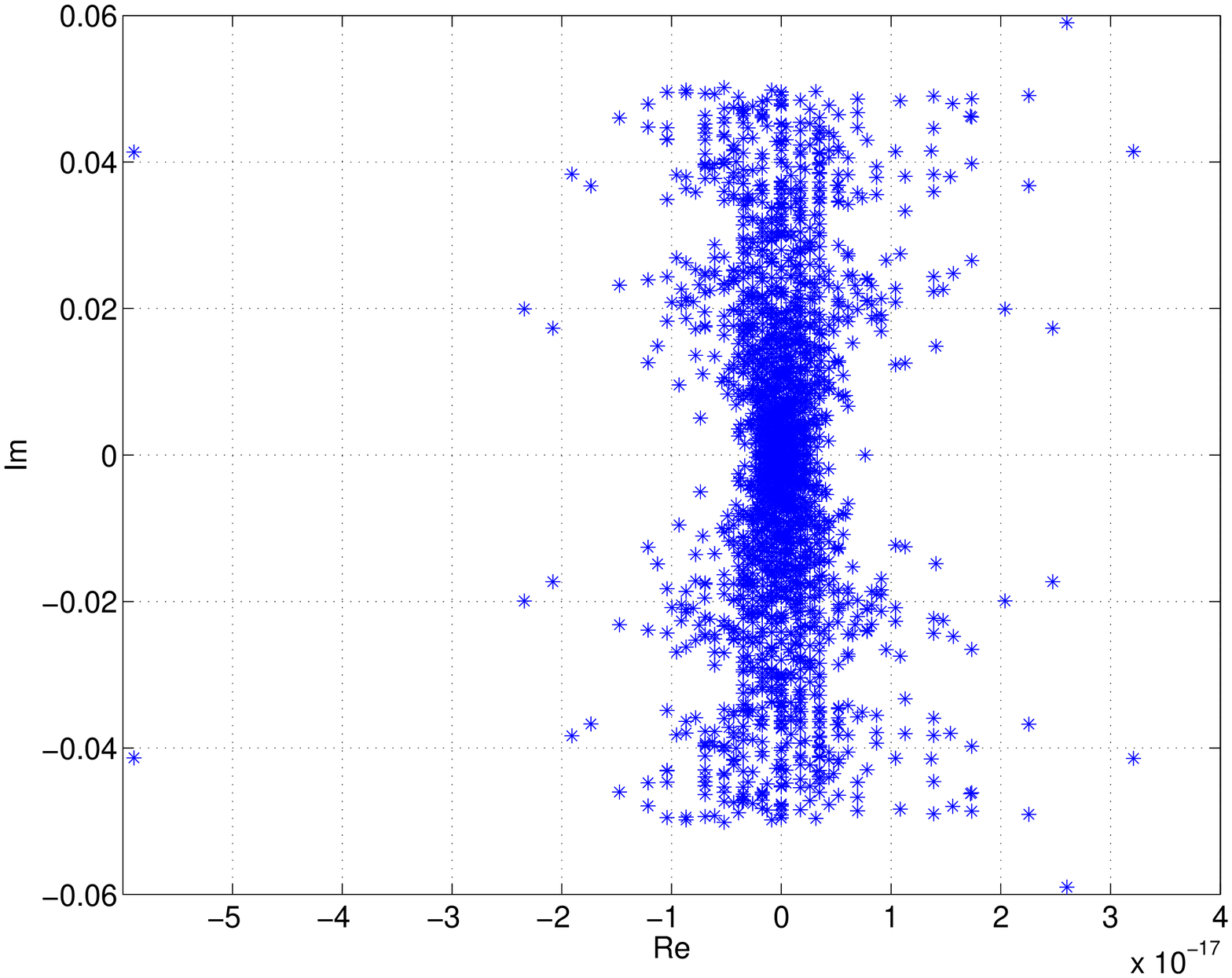}}
  \subfigure{\includegraphics[height=6.5cm,angle=-90]{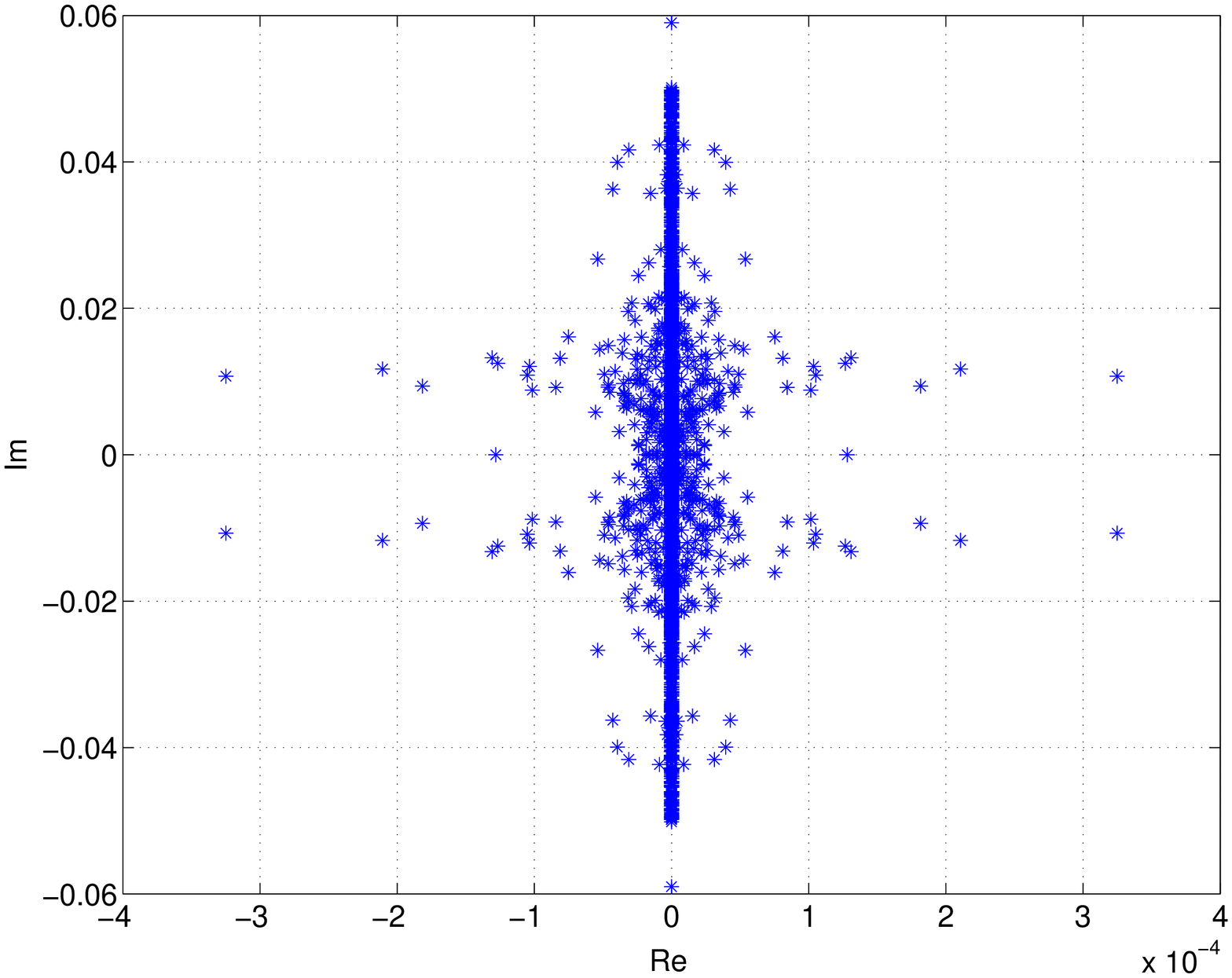}}
 \caption{The spectrum of the discrete approximation to (\ref{eq:PDE_S2}) multiplied by $h_R$ for the fourth-order discretization, with  stable SBP-preserving  interpolation (left) and unstable non-SBP interpolation (right). $15^2$ unknowns in the left (coarse) block. (Notice the different scales.)}\label{fig:spectra_4}
\end{figure}

\begin{figure}
  \subfigure{\includegraphics[height=6.5cm,angle=-90]{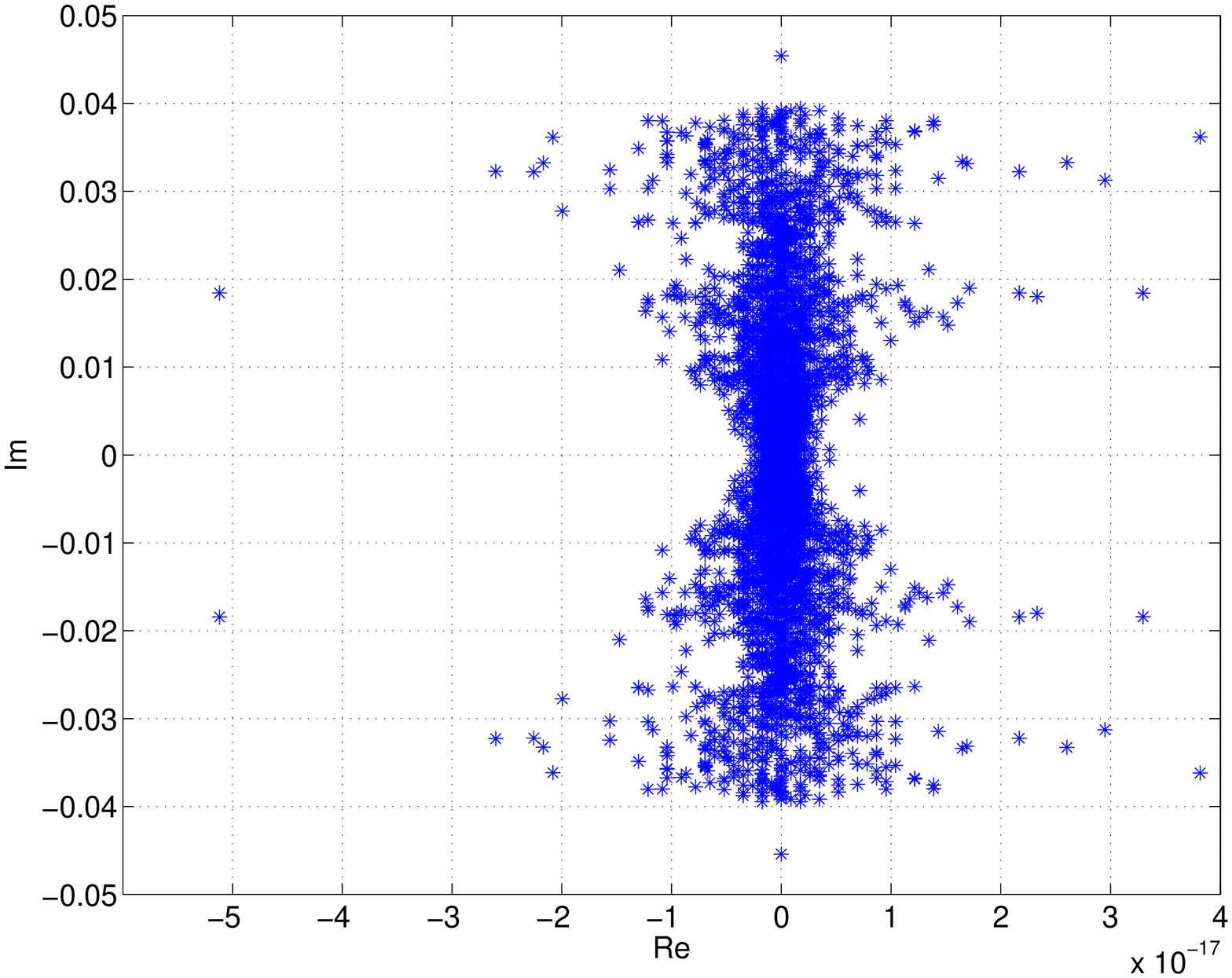}}
  \subfigure{\includegraphics[height=6.5cm,angle=-90]{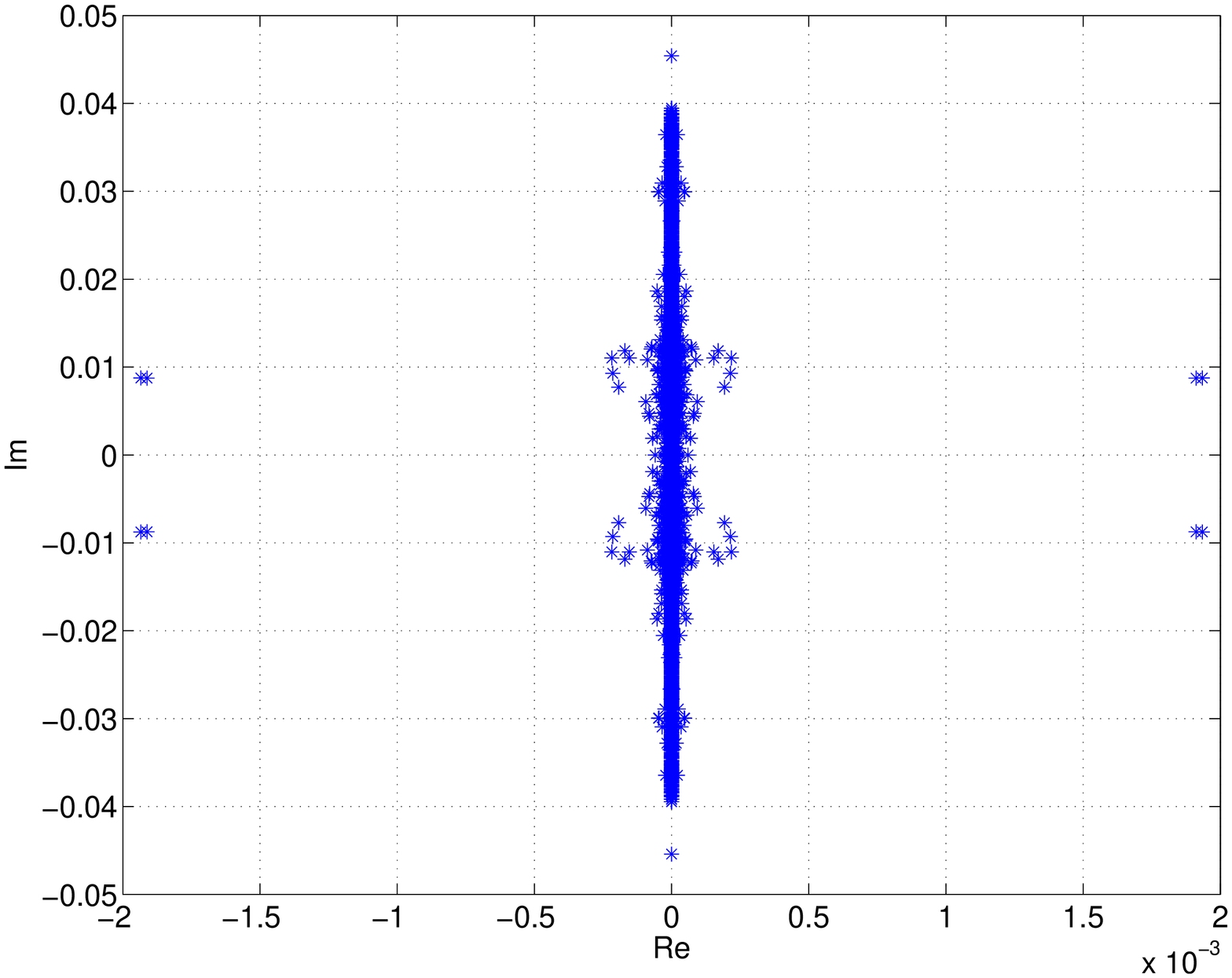}}
 \caption{The spectrum of the discrete approximation to (\ref{eq:PDE_S2}) multiplied by $h_R$ for the sixth-order discretization, with  stable SBP interpolation (left) and unstable non-SBP interpolation (right). $ 21^2$ unknowns in the left (coarse) block. (Notice the different scales.)}\label{fig:spectra_6}
\end{figure}

The eigenvalues of the numerical approximation of  the continuous test-problem (given  by (\ref{eq:hyper_pde}) and (\ref{eq:PDE_S2})) employing SBP-preserving interpolation operators fulfilling the conditions in (\ref{inter1}) for the fourth (using $15^2$ unknowns in the left (coarse mesh) region) and the sixth (using $21^2$ unknowns in the left (coarse mesh) region) order accurate cases are shown in Figures \ref{fig:spectra_4} and  \ref{fig:spectra_6} respectively. We have multiplied the eigenvalues with the grid-size ($h_R$) of the fine-mesh domain.  (The interpolation operators are presented in the Appendix.) To further indicate the importance of Lemma \ref{lemma:lemma1} we also present the corresponding results employing interpolation operators that do not fulfill the conditions in (\ref{inter1}) (and hence a violation to Lemma \ref{lemma:lemma1}), referred to as non-SBP interpolation. 

The non-SBP interpolation operators as compared to the corresponding SBP-preserving interpolation operator fulfilling the conditions in (\ref{inter1}) differ only at the boundaries. The boundary closures of the non-SBP interpolation operators (presented in the Appendix) are defined by the following two requirement: 1) to have the same formal accuracy at all grid-points, and 2) to use a minimal stencil width. 
The non-SBP interpolation operators clearly introduce an instability due to eigenvalues to the right of the imaginary axis. 

\subsection{Vortex in Free Space}\label{sec:freespace}
\begin{figure}
  \subfigure[t=0]{\includegraphics[height=12.5cm,angle=-90]{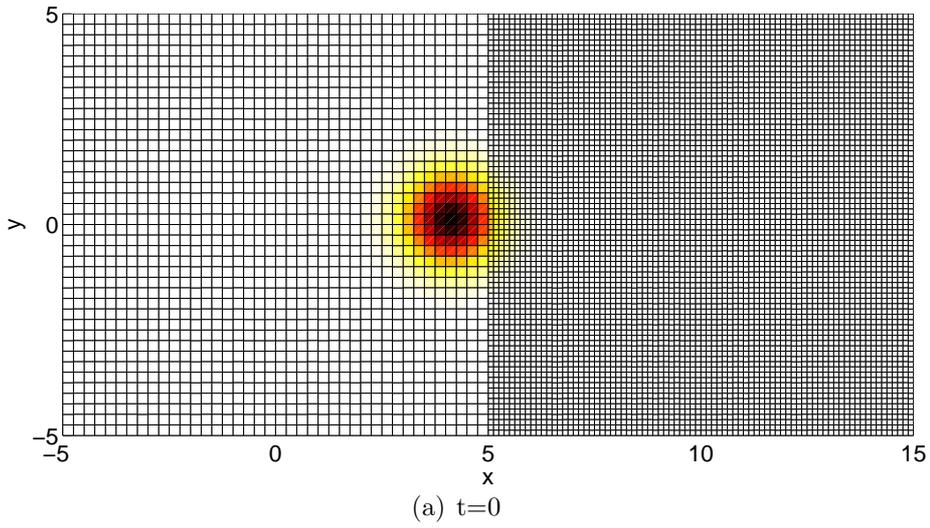}}
  \subfigure[t=1]{\includegraphics[height=12.5cm,angle=-90]{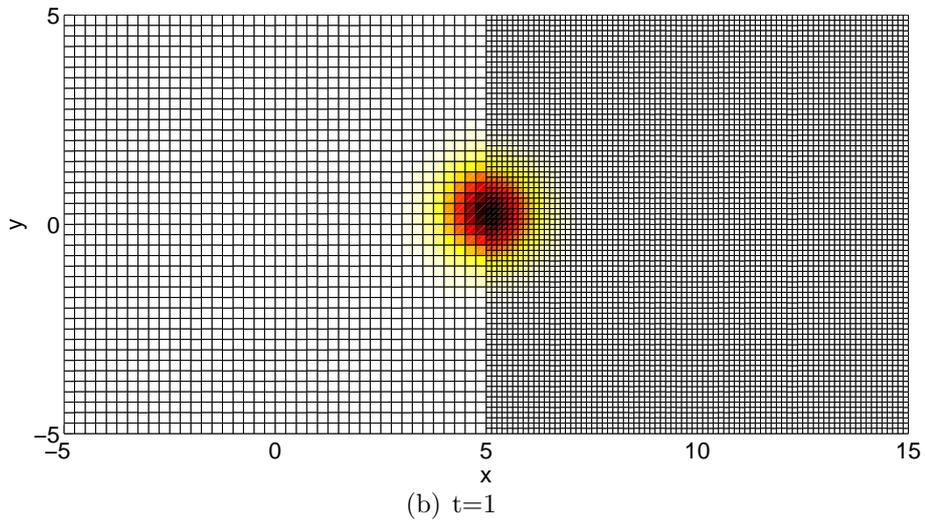}}
 \caption{Pressure contours for the vortex at t=0 (initial data)  and at $ t = 1$, computed with the sixth-order method using $41 \times 41$ grid points in the coarse (left domain).}\label{fig:solution}
\end{figure}

To test the accuracy of the present method a 2-D Euler-vortex (an analytic  solution \cite{MattssonSvard06} to the compressible Euler equations) is run across a multi-block interface. The problem consists of 2 blocks ($5\times5$ unit area)  having non-matching gridlines. The blocks are patched together according to (\ref{eq:hyper_SAT}). We set the Mach number $Ma=0.3$. The vortex is initiated at x=4 (one unit to the left of the interface at x=5) with a 10 degree angle of the background free-stream and then propagated to $t=1$, see Figure \ref{fig:solution}. 

The convergence rate is calculated as
\begin{equation}\label{eq:noggr}
q=\log_{10}\left(\frac{\|w-w^{(h_1)}\|_h}{\|w-w^{(h_2)}\|_h}\right)/\log_{10}\left(\frac{h_1}{h_2}\right)\;,
\end{equation}
where $w$ is the analytic solution and $w^{(h_1)}$ the corresponding
numerical solution with grid size $h_1$. $\|w-w^{(h_1)}\|_h$ is the discrete $l_2$ norm of the error. Since the solution consists of 2 blocks we add together the two different  $l_2$-errors. The standard explicit  4th-order Runge-Kutta method is  used for time integration.

The convergence results going from coarse to fine grid  using second, fourth, sixth and eighth-order accurate schemes are presented in Tables \ref{table:convergence_Left_second}-\ref{table:convergence_Left_eighth}. (The number $M$ is the number of unknowns in the y-direction in the coarse grid.)
\begin{table}
\centering
\begin{tabular}{ccccccccc}
\hline
$M$ &  $log(\rho)$ &$q^{(\rho)} $ & $log(l_2(u))$& $q^{(u)}$& $log(l_2(v))$ &$q^{(v)}$ & $log(l_2(p))$ &$q^{(p)} $\\
\hline
$ 51$ & -4.16 &  --  & -3.43 &  --  & -3.16 &  --  & -3.22 &  --  \\
$101$ & -4.76 & 2.00 & -4.04 & 2.02 & -3.78 & 2.04 & -3.82 & 2.01 \\
$151$ & -5.12 & 2.00 & -4.39 & 2.01 & -4.13 & 2.02 & -4.18 & 2.00 \\
\hline
\end{tabular}
\caption{$l_2$-error and convergence rate q. Second order case. Coarse to fine grid. $Ma=0.3$. Vortex centered at interface at $t=1$.}
\label{table:convergence_Left_second}
\end{table}

\begin{table}
\centering
\begin{tabular}{ccccccccc}
\hline
$M$ &  $log(\rho)$ &$q^{(\rho)} $ & $log(l_2(u))$& $q^{(u)}$& $log(l_2(v))$ &$q^{(v)}$ & $log(l_2(p))$ &$q^{(p)} $\\
\hline
$ 51$ & -4.61 &  --  & -4.22 &  --  & -3.85 &  --  & -3.80 &  --  \\
$101$ & -5.53 & 3.07 & -5.12 & 2.99 & -4.79 & 3.13 & -4.70 & 2.99 \\
$151$ & -6.07 & 3.05 & -5.65 & 2.99 & -5.33 & 3.08 & -5.23 & 3.01 \\
\hline
\end{tabular}
\caption{$l_2$-error and convergence rate q. Fourth order case. Coarse to fine grid. $Ma=0.3$. Vortex centered at interface at $t=1$.}
\label{table:convergence_Left_fourth}
\end{table}

\begin{table}
\centering
\begin{tabular}{ccccccccc}
\hline
$M$ &  $log(\rho)$ &$q^{(\rho)} $ & $log(l_2(u))$& $q^{(u)}$& $log(l_2(v))$ &$q^{(v)}$ & $log(l_2(p))$ &$q^{(p)} $\\
\hline
$ 51$ & -5.13 &  --  & -4.38 &  --  & -4.04 &  --  & -4.25 &  --  \\
$101$ & -6.27 & 3.79 & -5.66 & 4.25 & -5.30 & 4.19 & -5.53 & 4.24 \\
$151$ & -7.02 & 4.24 & -6.44 & 4.44 & -6.09 & 4.48 & -6.29 & 4.33 \\
\hline
\end{tabular}
\caption{$l_2$-error and convergence rate q. Sixth order case. Coarse to fine grid. $Ma=0.3$. Vortex centered at interface at $t=1$.}
\label{table:convergence_Left_sixth}
\end{table}

\begin{table}
\centering
\begin{tabular}{ccccccccc}
\hline
$M$ &  $log(\rho)$ &$q^{(\rho)} $ & $log(l_2(u))$& $q^{(u)}$& $log(l_2(v))$ &$q^{(v)}$ & $log(l_2(p))$ &$q^{(p)} $\\
\hline
$ 51$ & -4.72 &  --  & -4.44 &  --  & -3.93 &  --  & -4.00 &  --  \\
$101$ & -6.07 & 4.48 & -5.84 & 4.66 & -5.39 & 4.83 & -5.30 & 4.30 \\
$151$ & -6.93 & 4.87 & -6.67 & 4.75 & -6.27 & 5.04 & -6.13 & 4.74 \\
\hline
\end{tabular}
\caption{$l_2$-error and convergence rate q. Eighth order case. Coarse to fine grid. $Ma=0.3$. Vortex centered at interface at $t=1$.}
\label{table:convergence_Left_eighth}
\end{table}

The potential gain of a hybrid approach (here referring to the multi-block coupling  of  a higher order scheme on a coarse-grid domain to a lower order scheme on a fine-grid domain) is found by comparing the results (which are almost identical) in Tables  \ref{table:convergence_eigth_fourth} and \ref{table:convergence_Left_eighth}. In Table  \ref{table:convergence_eigth_fourth} we couple an eighth-order accurate SBP scheme in the left (coarse) domain to a fourth-order accurate SBP scheme in the right (fine) domain. In Table \ref{table:convergence_Left_eighth} the eighth-order stencil is used in both domains. 
In Table \ref{table:convergence_fourth_second} we couple a fourth-order accurate SBP scheme in the left (coarse) domain to a second-order accurate SBP scheme in the right (fine) domain. A comparison of  Tables \ref{table:convergence_fourth_second} and \ref{table:convergence_Left_fourth} indicates the benefit of employing  higher (than second) order accurate methods for wave dominated problems (the error of the second order schemes clearly dominates). 


\begin{table}
\centering
\begin{tabular}{ccccccccc}
\hline
$M$ &  $log(\rho)$ &$q^{(\rho)} $ & $log(l_2(u))$& $q^{(u)}$& $log(l_2(v))$ &$q^{(v)}$ & $log(l_2(p))$ &$q^{(p)} $\\
\hline
$ 51$ & -4.53 &  --  & -4.12 &  --  & -3.75 &  --  & -3.72 &  --  \\
$101$ & -5.28 & 2.51 & -4.81 & 2.30 & -4.47 & 2.39 & -4.45 & 2.45 \\
$151$ & -5.68 & 2.23 & -5.18 & 2.11 & -4.85 & 2.13 & -4.84 & 2.21 \\
\hline
\end{tabular}
\caption{$l_2$-error and convergence rate q. Fourth to second-order case. Coarse to fine grid. $Ma=0.3$. Vortex centered at interface at $t=1$.}
\label{table:convergence_fourth_second}
\end{table}

\begin{table}
\centering
\begin{tabular}{ccccccccc}
\hline
$M$ &  $log(\rho)$ &$q^{(\rho)} $ & $log(l_2(u))$& $q^{(u)}$& $log(l_2(v))$ &$q^{(v)}$ & $log(l_2(p))$ &$q^{(p)} $\\
\hline
$ 51$ & -4.74 &  --  & -4.43 &  --  & -3.93 &  --  & -4.01 &  --  \\
$101$ & -6.07 & 4.41 & -5.82 & 4.60 & -5.38 & 4.81 & -5.29 & 4.25 \\
$151$ & -6.89 & 4.69 & -6.63 & 4.63 & -6.25 & 4.96 & -6.11 & 4.65 \\
\hline
\end{tabular}
\caption{$l_2$-error and convergence rate q. Eighth to Fourth order case. Coarse to fine grid. $Ma=0.3$. Vortex centered at interface at $t=1$.}
\label{table:convergence_eigth_fourth}
\end{table}

\section{Conclusions}\label{sec:Conclusions}
Our approach have been  to use SBP operators, the SAT technique and SBP-preserving interpolation operators to enforce the interface conditions in a stable and accurate way for  general hyperbolic (and parabolic) problems, defined on nonconforming grids. 
The main objective was to construct interpolation operators that combined the following desirable properties:
\begin{itemize}
\item Stability by construction without interfering with the existing SBP schemes.
\item Maintaining  the overall convergence rate.
\item Maintaining simplicity of the numerical scheme.
\end{itemize}
To achieve the three properties above, we have constructed interpolation operators with the following
requirements: i) They preserve the summation by parts rule based on
the underlying  SBP scheme.  ii) They have the same order of accuracy as the underlying  SBP scheme.  iii) They are of minimal width in the interior.

Numerical computations for the 2-D compressible Euler equations  corroborate the stability and accuracy properties and also show that  a careful boundary treatment is necessary to guarantee stability of the interface coupling.

\section{Acknowledgments}\label{sec:Acknowledgments}

We would like to thank  Jan Nordstr\"om for very useful discussions, comments and suggestions throughout this work. Jan Nordstr\"om have made significant contributions in the framework of hybrid multi-block couplings of various SBP schemes utilizing the SAT technique. His earlier contributions have been a vital part in the development of the present study.  

\appendix

\section{Numerical Validation of Lemma \ref{lemma:lemma2} and Lemma \ref{lemma:lemma3}}\label{sec:energy2}
In this section we verify numerically the necessity of the extra conditions (\ref{inter2}) that turned up in Lemma \ref{lemma:lemma2}. We will also verify the modified penalties (\ref{eq:hyper_stability2}) that removed the necessity of  conditions (\ref{inter2}), in order to obtain an energy estimate, yet introducing a dissipative coupling (referring to the energy estimate (\ref{eq:norms2})). 

Consider (\ref{eq:hyper_SAT}) with the setup given by (\ref{eq:PDE_S2}). 
The scheme (\ref{eq:hyper_SAT}) with the penalty terms given by (\ref{eq:hyper_stability}) was proven stable if  $ \Sigma_{L}=A^- \;,\;  \Sigma_R=-A^+ $, assuming that both the conditions in (\ref{inter1}) and (\ref{inter2}) hold. However, for the sixth-order case condition (\ref{inter2}) could not be met. This led us to consider the scheme (\ref{eq:hyper_SAT}) with the penalty terms instead given by (\ref{eq:hyper_stability2}), where stability follows if $\Omega$ is symmetric and positive definite. In the present study we use $\Omega = \bar{A}$. The two different schemes (that differ only in how the SAT coupling is done) will be referred to as the characteristic coupling  and the quadratic coupling.

\begin{figure}
  \subfigure{\includegraphics[height=6.5cm,angle=-90]{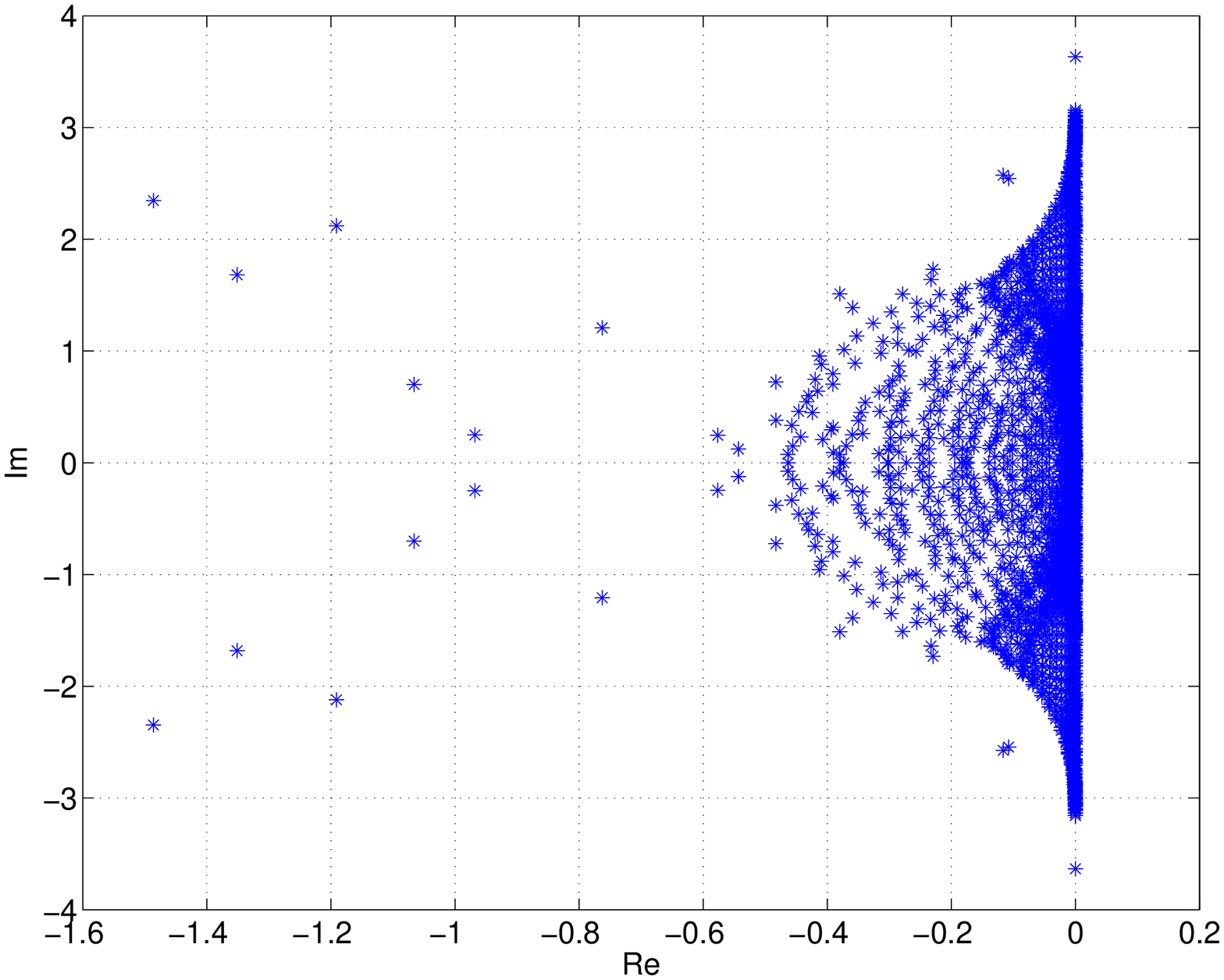}}
  \subfigure{\includegraphics[height=6.5cm,angle=-90]{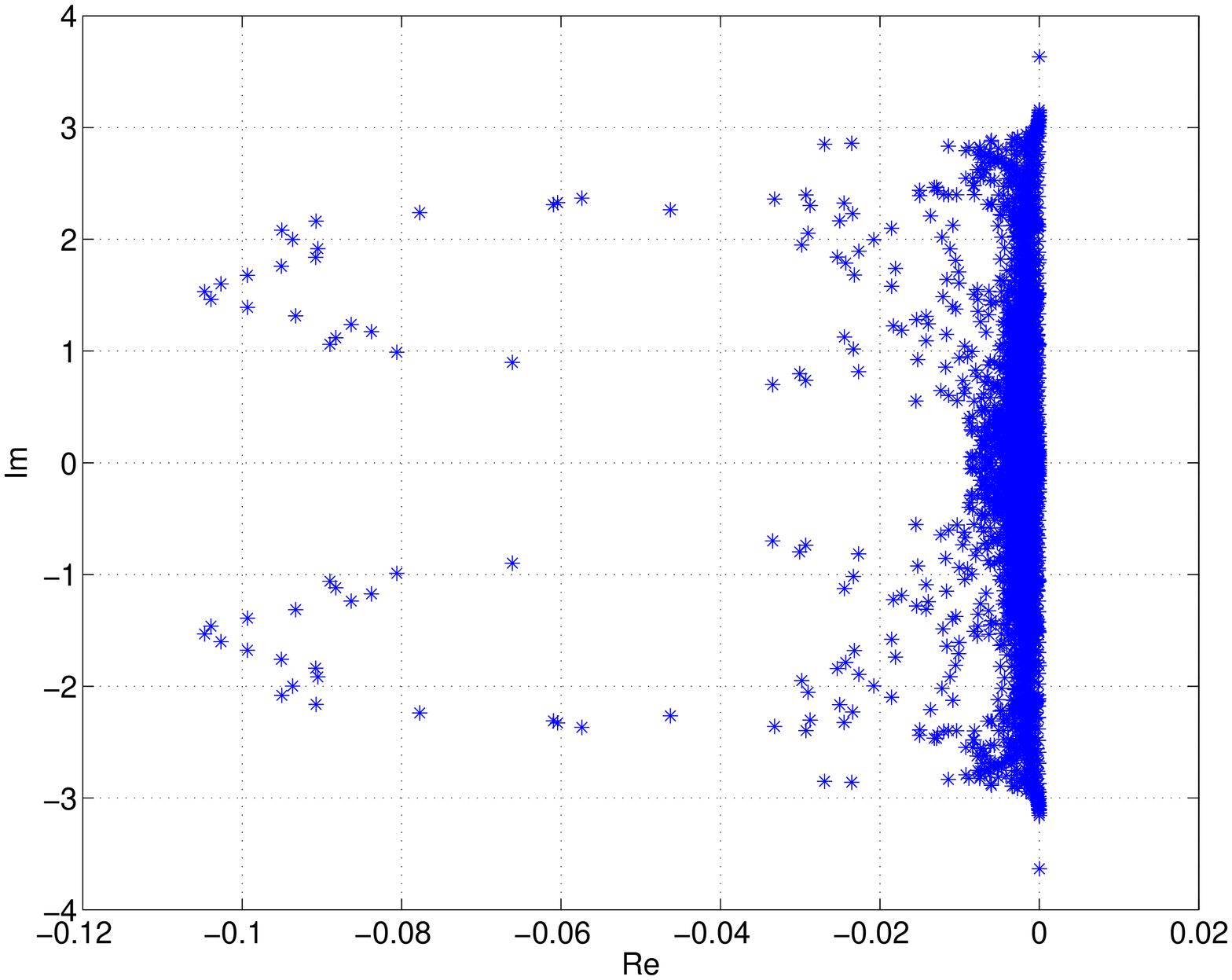}}
 \caption{The spectrum of the discrete approximation to (\ref{eq:PDE_S2}) multiplied by $h_R$ for the sixth-order discretization, with the characteristic coupling (left) and quadratic coupling (right). $21^2$ unknowns in the left (coarse) block.)}\label{fig:spectra_6a}
\end{figure}

The eigenvalues of the two different numerical approximations (the characteristic coupling  and the quadratic coupling) employing the SBP-preserving interpolation operators fulfilling the conditions in (\ref{inter1}) for  the sixth- (using $21^2$ unknowns in the left (coarse mesh) region) order accurate cases are shown in Figure \ref{fig:spectra_6a}  (compare with the corresponding non-dissipative coupling, presented in Figure \ref{fig:spectra_6}). We have multiplied the eigenvalues with the grid-size ($h_R$) of the fine-mesh domain. A closer look at the eigenvalues for the characteristic coupling reveals some eigenvalues with a positive real part, the largest of order $10^{-4}$. This is a strong indication that when employing a characteristic coupling the conditions in (\ref{inter2}) is necesary. (The interpolation operators for the second and fourth order accurate cases do obey the extra conditions in (\ref{inter2}). We have verified, although not presented here, that the second- and fourth-order accurate cases does not have eigenvalues with positive real part for the characteristic coupling). The quadratic coupling (not required to obey (\ref{inter2}))  on the other hand results in eigenvalues having non-positive real parts, in agreement with Lemma \ref{lemma:lemma3}.

\section{Parabolic problems}\label{sec:parabolic}
In this section we extend the analysis to include parabolic  problems (such as Navier-Stokes equations) of the form
\begin{equation}\label{eq:hyper_pde2}
\begin{array}{l}
\medskip {\bf u}_t+A{\bf u}_x+B{\bf u}_y=C_{11}{\bf u}_{xx}+C_{12}{\bf u}_{xy}+C_{21}{\bf u}_{yx}+ C_{22}{\bf u}_{yy},\; {\bf x} \in \Omega^{l,\,i}_{d,\,u},\\
\medskip {\bf v}_t+A{\bf v}_x+B{\bf v}_y=C_{11}{\bf v}_{xx}+C_{12}{\bf v}_{xy}+C_{21}{\bf v}_{yx}+ C_{22}{\bf v}_{yy},\; {\bf x} \in \Omega^{i,\,r}_{d,\,u},\\
\end{array}\;,
\end{equation}
where $A$ and $B$ are symmetric $k\times k$ matrices. The continuity conditions at the interface are given by: $u=v, u_x=v_x, u_y=v_y$.   Parabolicity requires that,
\begin{equation}\label{eq:dissipation}
x^T(C+C^T)x \ge 0\;,
\end{equation}
where 
\[
C=
\begin{bmatrix}
C_{11}&C_{12}\\
C_{21}&C_{22}\\
\end{bmatrix}.
\]
A corresponding semi-discrete approximation of (\ref{eq:hyper_pde2}) can be written
\begin{equation}\label{eq:hyper_SAT2}
\begin{array}{l}
\medskip u_t+A \otimes {\bf D_L^x}u+B \otimes {\bf D_L^y}u=P_L +SAT^I_{IL}+SAT^V_{IL}+SAT^V_L\\
\medskip  v_t+A \otimes {\bf D^R_x}v+B \otimes  {\bf D^R_y}v=P_R+SAT^I_{IR}+SAT^V_{IR}+SAT^V_R\;,\\
\end{array}
\end{equation}
where
\[
\begin{array}{l}
\medskip P_L=C_{11} \otimes {\bf D_L^x}{\bf D_L^x}u+C_{12} \otimes {\bf D_L^x}{\bf D_L^y}u+C_{21} \otimes {\bf D_L^y}{\bf D_L^x}u+C_{22} \otimes {\bf D_L^y}{\bf D_L^y}u\\
\medskip P_R=C_{11} \otimes {\bf D_R^x}{\bf D_R^x}v+C_{12} \otimes {\bf D_R^x}{\bf D_R^y}v+C_{21} \otimes {\bf D_R^y}{\bf D_R^x}v+C_{22} \otimes {\bf D_R^y}{\bf D_R^y}v\;,\\
\end{array}
\]
\[
\begin{split}
\medskip SAT^V_{IL}=&\Sigma^2_LC_{11} \otimes ({\bf H^{x}_L)^{-1}} ({\bf D_L^x})^T{\bf e^x_N}\left\{u_N-I_k\otimes I_{F2C}v_0\right\}\\
\medskip &C_{11}\Sigma^3_L \otimes ({\bf H^{x}_L)^{-1}}{\bf e^x_N}\left\{(I_k\otimes {\bf D_L^x} u)_N-I_k\otimes I_{F2C}(I_k\otimes {\bf D_R^x} v)_0\right\}\\
\medskip &\Sigma^4_LC_{21} \otimes ({\bf H^{x}_L)^{-1}} ({\bf D_L^y})^T{\bf e^x_N}\left\{u_N-I_k\otimes I_{F2C}v_0\right\}\\
\medskip &C_{21}\Sigma^5_L \otimes ({\bf H^{x}_L)^{-1}}{\bf e^x_N}\left\{(I_k\otimes {\bf D_L^y} u)_N-I_k\otimes I_{F2C}(I_k\otimes {\bf D_R^y} v)_0\right\}\;,\\
\end{split}
\]
\[
\begin{split}
\medskip SAT^V_{IR}=&\Sigma^2_RC_{11} \otimes ({\bf H^{x}_R})^{-1}({\bf D_R^x})^T{\bf e^x_0}\left\{v_0-I_k\otimes I_{C2F}u_N\right\}\\
\medskip &C_{11}\Sigma^3_R \otimes ({\bf H^{x}_R)^{-1}}{\bf e^x_N}\left\{(I_k\otimes {\bf D_R^x} v)_0-I_k\otimes I_{C2F}(I_k\otimes {\bf D_L^x} u)_N\right\}\\
\medskip &\Sigma^4_RC_{21}\ \otimes ({\bf H^{x}_R})^{-1}({\bf D_R^y})^T{\bf e^x_0}\left\{v_0-I_k\otimes I_{C2F}u_N\right\}\\
\medskip &C_{21}\Sigma^5_R \otimes ({\bf H^{x}_R)^{-1}}{\bf e^x_N}\left\{(I_k\otimes {\bf D_R^y} v)_0-I_k\otimes I_{C2F}(I_k\otimes {\bf D_L^y} u)_N\right\}\;,\\
\end{split}
\]
$SAT^I_{IL}$ and $SAT^I_{IR}$ are given by Eq. (\ref{eq:hyper_stability}), or Eq. (\ref{eq:hyper_stability2})).
$SAT^V_{IL}$ and $SAT^V_{IR}$ are the SAT penalty terms handling the viscous multi-block coupling. $SAT^V_{L,\,R}$ denote the viscous penalty terms at the outer boundaries. (These are tuned to obtain stability, see \cite{NordstromCarpenter01,MattssonSvard07,Svard04_Mattsson} for details.)

\begin{lemma}\label{lemma:lemma4}
The scheme (\ref{eq:hyper_SAT2}) is stable if  
\[
\begin{array}{l}
\medskip \Sigma^2_R= \Sigma^2 \;,\;  \Sigma^2_L=\Sigma^2-I  \;,\;  \Sigma^3_L=-\Sigma^2  \;,\;  \Sigma^3_R=I-\Sigma^2  \;,\\
\medskip \Sigma^4_R= \Sigma^4 \;,\;  \Sigma^4_L=\Sigma^4-I  \;,\;  \Sigma^5_L=-\Sigma^4  \;,\;  \Sigma^5_R=I-\Sigma^4  \;,\\
\end{array}\
\]
 and Lemma \ref{lemma:lemma1}, Lemma \ref{lemma:lemma2} or Lemma \ref{lemma:lemma3} hold, assuming that the boundary terms $SAT^V_{L,\,R}$ are correctly implemented.  $\Sigma^{2,\,4}$ are symmetric matrices.
\end{lemma} 

\begin{proof}
Apply the energy method by multiplying the first and second equation in (\ref{eq:hyper_SAT2}) by $u^TM_L$ and $v^TM_R$  (defined by (\ref{eq:norms}))  respectively, yielding,
\[
 \frac{d}{dt}(\|u\|_{M_L}+\|v\|_{M_R})=BT^V+IT^I+IT^V+(IT^V)^T\,.
\]
$BT^V$ corresponds to the outer boundary terms (which we assume are correctly implemented, i.e., bounded). $IT^I$ is given by (\ref{eq:hyper_terms}) and corresponds to the inviscid part. $IT^V+(IT^V)^T$ correspond to the viscous part  of the multi-block coupling given by:
\[
IT^V=
\begin{array}{l}
\medskip +u^T_NC_{11}(I_k+\Sigma^2_L+\Sigma^3_L)(H_L^y\otimes {\bf D_L^x} u)_N\\
\medskip -u^T_NC_{11}\otimes I^y(\Sigma^3_L\otimes H_L^y+\Sigma^2_R\otimes I_{C2F}^TH_R^y)(I_k\otimes {\bf D_R^x} v)_0\\
\medskip -v^T_0C_{11}(I_k-\Sigma^2_R-\Sigma^3_R)(H_R^y\otimes {\bf D_R^x} v)_0\\
\medskip -v^T_0C_{11}\otimes I^y(\Sigma^3_R\otimes H_R^y+\Sigma^2_L\otimes I_{F2C}^TH_L^y)(I_k\otimes {\bf D_L^x} u)_N\\
\medskip +u^T_NC_{21}(I_k+\Sigma^4_L+\Sigma^5_L)(H_L^y\otimes {\bf D_L^y} u)_N\\
\medskip -u^T_NC_{21}\otimes I^y(\Sigma^5_L\otimes H_L^y+\Sigma^4_R\otimes I_{C2F}^TH_R^y)(I_k\otimes {\bf D_R^y} v)_0\\
\medskip -v^T_0C_{21}(I_k-\Sigma^4_R-\Sigma^5_R)(H_R^y\otimes {\bf D_R^y} v)_0\\
\medskip -v^T_0C_{21}\otimes I^y(\Sigma^5_R\otimes H_R^y+\Sigma^4_L\otimes I_{F2C}^TH_L^y)(I_k\otimes {\bf D_L^y} u)_N\;.\\
\end{array}\;,
\]
A  stable coupling is obtained if $IT^I$ is non-positive and if $IT^V$ vanish. 
By assumption Lemma \ref{lemma:lemma1}, Lemma \ref{lemma:lemma2} or Lemma \ref{lemma:lemma3} hold which means that  $IT^I$ is non-positive. By choosing:
\[
\begin{array}{l}
\medskip \Sigma^2_R= \Sigma^2 \;,\;  \Sigma^2_L=\Sigma^2-I_k  \;,\;  \Sigma^3_L=-\Sigma^2  \;,\;  \Sigma^3_R=I_k-\Sigma^2  \;,\\
\medskip \Sigma^4_R= \Sigma^4 \;,\;  \Sigma^4_L=\Sigma^4-I_k  \;,\;  \Sigma^5_L=-\Sigma^4  \;,\;  \Sigma^5_R=I_k-\Sigma^4  \;,\\
\end{array}\
\]
(in combination with (\ref{inter1})) $IT^V=0$.
 $\square$
\end{proof}

\begin{remark}
The symmetric matrices $\Sigma^2$ and $\Sigma^4$ can be chosen arbitrarily for stability. (Note that $\Sigma^{2\,,4}$ are not powers of $\Sigma$, but different parameters). However, the specific choice will affect the eigenvalues of  (\ref{eq:hyper_SAT2}). We use $\Sigma^2=\Sigma^4=I_k/2$ since it leads to a compact spectrum (we have no proof that this particular choice is optimal).
\end {remark}

\section{Interpolation Operators}

\subsection{Second-Order case}
The discrete norm is given by $H=h\,diag(\frac{1}{2},1,\ldots,1,\frac{1}{2})$.  
The interpolation operators are given by:
\[
I_{F2C}=
\left[ \begin {array}{cccccccccc}
 {\frac {11}{20}}&\frac{1}{2}&-\frac{1}{
20}&&&&&&&\\
\noalign{\medskip}-\frac{1}{40}&\frac{1}{4}&{\frac {11}{20}}&\frac{1}{4}&-\frac{1}{40}&&&&&\\
\noalign{\medskip}&&-\frac{1}{40}&\frac{1}{4}&{\frac {11}{20}}&\frac{1}{4}&-\frac{1}{40}&&\\
\noalign{\medskip}&&&\ddots&\ddots&\ddots&\ddots&\ddots&&\\
\noalign{\medskip}&&&-\frac{1}{40}&\frac{1}{4}&{\frac {11}{20}}&\frac{1}{4}&-\frac{1}{40}&&\\
\noalign{\medskip}&&&&&-\frac{1}{40}&\frac{1}{4}&{\frac {11}{20}}&\frac{1}{4}&-\frac{1}{40}\\
\noalign{\medskip}&&&&&&&-\frac{1}{20}&\frac{1}{2}&{\frac {11}{20}}\end {array} \right] 
\]
and
\[
I_{C2F}=
\left[ \begin {array}{ccccccc} {\frac {11}{10}}&-\frac{1}{10}&&&&&\\
\noalign{\medskip}\frac{1}{2}&\frac{1}{2}&&&&&\\
\noalign{\medskip}-\frac{1}{20}&{\frac {11}{10}}&-\frac{1}{20}&&&&\\
\noalign{\medskip}&\frac{1}{2}&\frac{1}{2}&&&&\\
\noalign{\medskip}&-\frac{1}{20}&{\frac {11}{10}}&-\frac{1}{20}&&&\\
\noalign{\medskip}&&\frac{1}{2}&\frac{1}{2}&&&\\
\noalign{\medskip}&&-\frac{1}{20}&{\frac {11}{10}}&-\frac{1}{20}&&\\
\noalign{\medskip}&&&\ddots&\ddots&\ddots&\\
\noalign{\medskip}&&&&-\frac{1}{20}&{\frac {11}{10}}&-\frac{1}{20}\\
\noalign{\medskip}&&&&&\frac{1}{2}&\frac{1}{2}\\
\noalign{\medskip}&&&&&-\frac{1}{10}&{\frac {11}{10}}
\end {array} \right] 
\]

\subsection{Fourth- to Second-Order Case}

Fourth to second-order, with second-order accurate interpolation. The upper left corner of $I_{F2C}$ is given by:

\[
{\scriptsize
\begin{array}{lll}
\begin{array}{l}
\medskip a_{1,1}= {\frac {198814379483797276626589}{292154366866608903386020}}\\
\medskip a_{1,2}= {\frac {33159311983113441321394}{73038591716652225846505}}\\
\medskip a_{1,3}= {\frac {1971623072622604595228923}{16828091531516672835034752}}\\
\medskip a_{1,4}= -{\frac {95808836247685639307767}{824906447623366315442880}}\\
\medskip a_{1,5}= -{\frac {1754070195140572218604399}{10517557207197920521896720}}\\
\medskip a_{1,6}= -{\frac {237072760963890638174137}{14023409609597227362528960}}\\
\medskip a_{1,7}= {\frac {55326079761477589117471}{2804681921919445472505792}}\\
\medskip a_{1,8}= -{\frac {120787724660864994449209}{14023409609597227362528960}}\\
\medskip a_{1,9}= {\frac {176943055421257964368859}{10517557207197920521896720}}\\
\medskip a_{1,10}= {\frac {32831043571234596789193}{14023409609597227362528960}}\\
\medskip a_{1,11}= {\frac {88481545475605230181351}{4949438685740197892657280}}\\
\end {array}
&
\begin{array}{l}
\medskip a_{2,1}={\frac {1342353899420590372183}{101394750853705442939854}}\\
\medskip a_{2,2}={\frac {87193543695192148191232}{253486877134263607349635}}\\
\medskip a_{2,3}={\frac {5578077558511726076198789}{19467792163911445044451968}}\\
\medskip a_{2,4}={\frac {3895790315853711272834919}{16223160136592870870376640}}\\
\medskip a_{2,5}={\frac {2244981324267232772622413}{12167370102444653152782480}}\\
\medskip a_{2,6}=-{\frac {213524706698808911938359}{16223160136592870870376640}}\\
\medskip a_{2,7}=-{\frac {2744378984383218872921257}{48669480409778612611129920}}\\
\medskip a_{2,8}={\frac {177648545875941908529241}{16223160136592870870376640}}\\
\medskip a_{2,9}=-{\frac {45572983468520710671107}{12167370102444653152782480}}\\
\medskip a_{2,0}={\frac {3876549717482752134187}{3244632027318574174075328}}\\
\medskip a_{2,11}=-{\frac {704697674731763751811223}{97338960819557225222259840}}\\
\end {array}
&
\begin{array}{l}
\medskip a_{3,1}=-{\frac {6011306572161601487699}{738978692662598990917580}}\\
\medskip a_{3,2}=-{\frac {17239549725384800082266}{184744673165649747729395}}\\
\medskip a_{3,3}={\frac {13200899856172905887766503}{70941954495609503128087680}}\\
\medskip a_{3,4}={\frac {3693069541319881367608217}{11823659082601583854681280}}\\
\medskip a_{3,5}={\frac {1607455065911582798562599}{8867744311951187891010960}}\\
\medskip a_{3,6}={\frac {3514316393052784877273751}{11823659082601583854681280}}\\
\medskip a_{3,7}={\frac {8869555969232379886041937}{35470977247804751564043840}}\\
\medskip a_{3,8}=-{\frac {652819712407377597489161}{11823659082601583854681280}}\\
\medskip a_{3,9}=-{\frac {129041283038038640595559}{1773548862390237578202192}}\\
\medskip a_{3,10}=-{\frac {21210275952319643393399}{11823659082601583854681280}}\\
\medskip a_{3,11}={\frac {300890344598183653971671}{70941954495609503128087680}}\\
\end {array}
\end {array}}
\]

The interior of $I_{F2C}$ is given by:

\[
\begin{array}{ll}
\begin{array}{l}
\medskip a_{4}=-{\frac {292764282684548619564311}{4949438685740197892657280}}\\
\medskip a_{3}=-{\frac {27348461476516948621801}{824906447623366315442880}}\\
\medskip a_{2}={\frac {138809444447592059191981}{618679835717524736582160}}\\
\end {array}
&
\begin{array}{l}
\medskip a_{1}={\frac {233575073382358527482521}{824906447623366315442880}}\\
\medskip a_{0}={\frac {139882799512953873064261}{824906447623366315442880}}\\
\end {array}
\end {array}
\]

The upper left corner of $I_{C2F}$ is given by:

\[
{\scriptsize
\begin{array}{lll}
\begin{array}{l}
\medskip b_{1,1}= {\frac {198814379483797276626589}{206226611905841578860720}}\\
\medskip b_{1,2}={\frac {1342353899420590372183}{20622661190584157886072}}\\
\medskip b_{1,3}=-{\frac {6011306572161601487699}{206226611905841578860720}}\\
\medskip b_{2,1}={\frac {16579655991556720660697}{51556652976460394715180}}\\
\medskip b_{2,2}={\frac {10899192961899018523904}{12889163244115098678795}}\\
\medskip b_{2,3}=-{\frac {8619774862692400041133}{51556652976460394715180}}\\
\medskip b_{3,1}={\frac {1971623072622604595228923}{23757305691552949884754944}}\\
\medskip b_{3,2}={\frac {5578077558511726076198789}{7919101897184316628251648}}\\
\medskip b_{3,3}={\frac {13200899856172905887766503}{39595509485921583141258240}}\\
\medskip b_{3,4}=-{\frac {14345449851542882358651239}{118786528457764749423774720}}\\
\medskip b_{4,1}=-{\frac {1628750216210655868232039}{19797754742960791570629120}}\\
\medskip b_{4,2}={\frac {1298596771951237090944973}{2199750526995643507847680}}\\
\medskip b_{4,3}={\frac {3693069541319881367608217}{6599251580986930523543040}}\\
\medskip b_{4,4}=-{\frac {1340074612349330482468249}{19797754742960791570629120}}\\
\medskip b_{5,1}=-{\frac {1754070195140572218604399}{14848316057220593677971840}}\\
\medskip b_{5,2}={\frac {2244981324267232772622413}{4949438685740197892657280}}\\
\medskip b_{5,3}={\frac {1607455065911582798562599}{4949438685740197892657280}}\\
\medskip b_{5,4}={\frac {6801662777932010900407069}{14848316057220593677971840}}\\
\medskip b_{5,5}=-{\frac {292764282684548619564311}{2474719342870098946328640}}\\
\end {array}
&
\begin{array}{l}
\medskip b_{6,1}=-{\frac {237072760963890638174137}{19797754742960791570629120}}\\
\medskip b_{6,2}=-{\frac {71174902232936303979453}{2199750526995643507847680}}\\
\medskip b_{6,3}={\frac {1171438797684261625757917}{2199750526995643507847680}}\\
\medskip b_{6,4}={\frac {11445178595735567846643529}{19797754742960791570629120}}\\
\medskip b_{6,5}=-{\frac {27348461476516948621801}{412453223811683157721440}}\\
\medskip b_{7,1}={\frac {55326079761477589117471}{3959550948592158314125824}}\\
\medskip b_{7,2}=-{\frac {2744378984383218872921257}{19797754742960791570629120}}\\
\medskip b_{7,3}={\frac {8869555969232379886041937}{19797754742960791570629120}}\\
\medskip b_{7,4}={\frac {6854257176134739780148789}{19797754742960791570629120}}\\
\medskip b_{7,5}={\frac {138809444447592059191981}{309339917858762368291080}}\\
\medskip b_{7,6}=-{\frac {292764282684548619564311}{2474719342870098946328640}}\\
\medskip b_{8,1}=-{\frac {120787724660864994449209}{19797754742960791570629120}}\\
\medskip b_{8,2}={\frac {177648545875941908529241}{6599251580986930523543040}}\\
\medskip b_{8,3}=-{\frac {652819712407377597489161}{6599251580986930523543040}}\\
\medskip b_{8,4}={\frac {11445178595735567846643529}{19797754742960791570629120}}\\
\medskip b_{8,5}={\frac {233575073382358527482521}{412453223811683157721440}}\\
\medskip b_{8,6}=-{\frac {27348461476516948621801}{412453223811683157721440}}\\
\medskip b_{9,1}={\frac {176943055421257964368859}{14848316057220593677971840}}\\
\medskip b_{9,2}=-{\frac {45572983468520710671107}{4949438685740197892657280}}\\
\medskip b_{9,3}=-{\frac {129041283038038640595559}{989887737148039578531456}}\\
\end {array}
&
\begin{array}{l}
\medskip b_{9,4}={\frac {6801662777932010900407069}{14848316057220593677971840}}\\
\medskip b_{9,5}={\frac {139882799512953873064261}{412453223811683157721440}}\\
\medskip b_{9,6}={\frac {138809444447592059191981}{309339917858762368291080}}\\
\medskip b_{9,7}=-{\frac {292764282684548619564311}{2474719342870098946328640}}\\
\medskip b_{10,1}={\frac {32831043571234596789193}{19797754742960791570629120}}\\
\medskip b_{10,2}={\frac {3876549717482752134187}{1319850316197386104708608}}\\
\medskip b_{10,3}=-{\frac {21210275952319643393399}{6599251580986930523543040}}\\
\medskip b_{10,4}=-{\frac {1340074612349330482468249}{19797754742960791570629120}}\\
\medskip b_{10,5}={\frac {233575073382358527482521}{412453223811683157721440}}\\
\medskip b_{10,6}={\frac {233575073382358527482521}{412453223811683157721440}}\\
\medskip b_{10,7}=-{\frac {27348461476516948621801}{412453223811683157721440}}\\
\medskip b_{11,1}={\frac {1504186273085288913082967}{118786528457764749423774720}}\\
\medskip b_{11,2}=-{\frac {704697674731763751811223}{39595509485921583141258240}}\\
\medskip b_{11,3}={\frac {300890344598183653971671}{39595509485921583141258240}}\\
\medskip b_{11,4}=-{\frac {14345449851542882358651239}{118786528457764749423774720}}\\
\medskip b_{11,5}={\frac {138809444447592059191981}{309339917858762368291080}}\\
\medskip b_{11,6}={\frac {139882799512953873064261}{412453223811683157721440}}\\
\medskip b_{11,7}={\frac {138809444447592059191981}{309339917858762368291080}}\\
\medskip b_{11,8}=-{\frac {292764282684548619564311}{2474719342870098946328640}}\\
\end {array}
\end {array}}
\]

\subsection{Fourth-Order Case}
The discrete norm is given by $H=h\,diag(\frac{17}{48}, \frac{59}{48}, \frac{43}{48}, \frac{49}{48},1,\ldots)$.  

The upper left corner of $I_{F2C}$ is given by:
\[
{\scriptsize
\begin{array}{lll}
\begin{array}{l}
\medskip a_{1,1}={\frac {37625155150141581259}{74136705821905530032}}\\
\medskip a_{1,2}={\frac {109494000531368265009}{157540499871549251318}}\\
\medskip a_{1,3}={\frac {598984862036951996355}{20165183983558304168704}}\\
\medskip a_{1,4}=-{\frac {259928664916563226325}{2520647997944788021088}}\\
\medskip a_{1,5}=-{\frac {256587334010286818055}{5041295995889576042176}}\\
\medskip a_{1,6}=-{\frac {113845943840605140771}{2520647997944788021088}}\\
\medskip a_{1,7}=-{\frac {537477213186037277591}{10082591991779152084352}}\\
\medskip a_{1,8}=-{\frac {89579548577713740579}{2520647997944788021088}}\\
\medskip a_{1,9}=-{\frac {33361750807060821603}{5041295995889576042176}}\\
\medskip a_{1,10}={\frac {45876275221828488395}{2520647997944788021088}}\\
\medskip a_{1,11}={\frac {891498675213598216803}{20165183983558304168704}}\\
\end {array}
&
\begin{array}{l}
\medskip a_{2,1}=-{\frac {9465638066209876131}{2187032821746213135944}}\\
\medskip a_{2,2}={\frac {6477310271684014677}{18534176455476382508}}\\
\medskip a_{2,3}={\frac {23789903736542318389255}{69985050295878820350208}}\\
\medskip a_{2,4}={\frac {2109162960890226252711}{8748131286984852543776}}\\
\medskip a_{2,5}={\frac {776008005836532793785}{17496262573969705087552}}\\
\medskip a_{2,6}={\frac {46983667240315552089}{8748131286984852543776}}\\
\medskip a_{2,7}={\frac {1087641244818460785021}{34992525147939410175104}}\\
\medskip a_{2,8}={\frac {220860994180249341801}{8748131286984852543776}}\\
\medskip a_{2,9}={\frac {77236835126760009513}{17496262573969705087552}}\\
\medskip a_{2,10}=-{\frac {96386094557526072417}{8748131286984852543776}}\\
\medskip a_{2,11}=-{\frac {1793510683939834799913}{69985050295878820350208}}\\
\end {array}
&
\begin{array}{l}
\medskip a_{3,1}={\frac {9465638066209876131}{3187878350341937791376}}\\
\medskip a_{3,2}=-{\frac {54391101655540111975}{796969587585484447844}}\\
\medskip a_{3,3}={\frac {1629291546275511894441}{51006053605471004662016}}\\
\medskip a_{3,4}={\frac {2010507950879293630905}{6375756700683875582752}}\\
\medskip a_{3,5}={\frac {6208709747109065353983}{12751513401367751165504}}\\
\medskip a_{3,6}={\frac {1804441319301200448807}{6375756700683875582752}}\\
\medskip a_{3,7}=-{\frac {58210841472168154533}{25503026802735502331008}}\\
\medskip a_{3,8}=-{\frac {395393460093074051961}{6375756700683875582752}}\\
\medskip a_{3,9}=-{\frac {180548419983997949901}{12751513401367751165504}}\\
\medskip a_{3,10}={\frac {55143363449566679649}{6375756700683875582752}}\\
\medskip a_{3,11}={\frac {21221519586950580219}{1186187293150488480512}}\\
\end {array}
\end {array}}
\]

The interior stencil is given by:
\[
\begin{array}{lll}
\begin{array}{l}
\medskip a_4=-{\frac {10513333512638366307}{1186187293150488480512}}\\
\medskip a_3=-\frac{1}{32}\\
\end {array}
&
\begin{array}{l}
\medskip a_2={\frac {10513333512638366307}{296546823287622120128}}\\
\end {array}
&
\begin{array}{l}
\medskip a_1={\frac {9}{32}}\\
\medskip a_0{\frac {265006822749707021207}{593093646575244240256}}\\
\end {array}
\end {array}
\]

The upper left part of $I_{C2F}$ is given by:
\[
{\scriptsize
\begin{array}{lll}
\begin{array}{l}
\medskip b_{1,1}={\frac {37625155150141581259}{37068352910952765016}}\\
\medskip b_{1,2}=-{\frac {556802239188816243}{18534176455476382508}}\\
\medskip b_{1,3}={\frac {556802239188816243}{37068352910952765016}}\\
\medskip b_{2,1}={\frac {1855830517480818051}{4633544113869095627}}\\
\medskip b_{2,2}={\frac {6477310271684014677}{9267088227738191254}}\\
\medskip b_{2,3}=-{\frac {921883078907459525}{9267088227738191254}}\\
\medskip b_{3,1}={\frac {598984862036951996355}{25503026802735502331008}}\\
\medskip b_{3,2}={\frac {23789903736542318389255}{25503026802735502331008}}\\
\medskip b_{3,3}={\frac {1629291546275511894441}{25503026802735502331008}}\\
\medskip b_{3,4}=-{\frac {515153342119279949043}{25503026802735502331008}}\\
\medskip b_{4,1}=-{\frac {5304666630950269925}{74136705821905530032}}\\
\medskip b_{4,2}={\frac {43044142058984209239}{74136705821905530032}}\\
\medskip b_{4,3}={\frac {41030774507740686345}{74136705821905530032}}\\
\medskip b_{4,4}=-\frac{1}{16}\\
\medskip b_{5,1}=-{\frac {85529111336762272685}{2372374586300976961024}}\\
\medskip b_{5,2}={\frac {258669335278844264595}{2372374586300976961024}}\\
\medskip b_{5,2}={\frac {2069569915703021784661}{2372374586300976961024}}\\
\medskip b_{5,2}={\frac {171717780706426649681}{2372374586300976961024}}\\
\medskip b_{5,2}=-{\frac {10513333512638366307}{593093646575244240256}}\\
\end {array}
&
\begin{array}{l}
\medskip b_{6,1}=-{\frac {37948647946868380257}{1186187293150488480512}}\\
\medskip b_{6,2}={\frac {15661222413438517363}{1186187293150488480512}}\\
\medskip b_{6,2}={\frac {601480439767066816269}{1186187293150488480512}}\\
\medskip b_{6,2}={\frac {147}{256}}\\
\medskip b_{6,2}=-\frac{1}{16}\\
\medskip b_{7,1}=-{\frac {537477213186037277591}{14234247517805861766144}}\\
\medskip b_{7,2}={\frac {362547081606153595007}{4744749172601953922048}}\\
\medskip b_{7,3}=-{\frac {19403613824056051511}{4744749172601953922048}}\\
\medskip b_{7,4}={\frac {12985334314735644039143}{14234247517805861766144}}\\
\medskip b_{7,5}={\frac {10513333512638366307}{148273411643811060064}}\\
\medskip b_{7,6}=-{\frac {10513333512638366307}{593093646575244240256}}\\
\medskip b_{8,1}=-{\frac {29859849525904580193}{1186187293150488480512}}\\
\medskip b_{8,2}={\frac {73620331393416447267}{1186187293150488480512}}\\
\medskip b_{8,3}=-{\frac {131797820031024683987}{1186187293150488480512}}\\
\medskip b_{8,4}={\frac {147}{256}}\\
\medskip b_{8,5}={\frac {9}{16}}\\
\medskip b_{8,6}=-\frac{1}{16}\\
\medskip b_{9,1}=-{\frac {11120583602353607201}{2372374586300976961024}}\\
\medskip b_{9,2}={\frac {25745611708920003171}{2372374586300976961024}}\\
\end {array}
&
\begin{array}{l}
\medskip b_{9,3}=-{\frac {60182806661332649967}{2372374586300976961024}}\\
\medskip b_{9,4}={\frac {171717780706426649681}{2372374586300976961024}}\\
\medskip b_{9,5}={\frac {265006822749707021207}{296546823287622120128}}\\
\medskip b_{9,6}={\frac {10513333512638366307}{148273411643811060064}}\\
\medskip b_{9,7}=-{\frac {10513333512638366307}{593093646575244240256}}\\
\medskip b_{10,1}={\frac {45876275221828488395}{3558561879451465441536}}\\
\medskip b_{10,2}=-{\frac {32128698185842024139}{1186187293150488480512}}\\
\medskip b_{10,3}={\frac {18381121149855559883}{1186187293150488480512}}\\
\medskip b_{10,4}=-{\frac {49}{768}}\\
\medskip b_{10,5}={\frac {9}{16}}\\
\medskip b_{10,6}={\frac {9}{16}}\\
\medskip b_{10,7}=-\frac{1}{16}\\
\medskip b_{11,1}={\frac {297166225071199405601}{9489498345203907844096}}\\
\medskip b_{11,2}=-{\frac {597836894646611599971}{9489498345203907844096}}\\
\medskip b_{11,3}={\frac {304175114079624983139}{9489498345203907844096}}\\
\medskip b_{11,4}=-{\frac {171717780706426649681}{9489498345203907844096}}\\
\medskip b_{11,5}={\frac {10513333512638366307}{148273411643811060064}}\\
\medskip b_{11,6}={\frac {265006822749707021207}{296546823287622120128}}\\
\medskip b_{11,7}={\frac {10513333512638366307}{148273411643811060064}}\\
\medskip b_{11,8}=-{\frac {10513333512638366307}{593093646575244240256}}\\
\end {array}
\end {array}}
\]

\subsection{Fourth-Order Case, Non-SBP Interpolation}

The upper left part of $I_{C2F}$:
\[
\left( \begin {array}{ccccccc} 1& & & & && \\
\noalign{\medskip}{\frac {5}{16}}&{\frac {15}{16}}&-{\frac {5}{16}}&\frac{1}{16}& && \\
\noalign{\medskip} &1& & & & &\\
\noalign{\medskip} 2a_3&2a_1&2a_1&2a_3&&&\\
\noalign{\medskip} 2a_4&2a_2&2a_0&2a_2&2a_4&&\\
\noalign{\medskip} & &2a_3&2a_1&2a_1&2a_3&\\
\noalign{\medskip} & &2a_4&2a_2&2a_0&2a_2&2a_4\\
\noalign{\medskip} &&&\ddots&\ddots&\ddots&\\
\end{array}\right)
\]

The lower left part of $I_{F2C}$:
\[
\left( \begin {array}{ccccccccccc} 1& & & & & & & & & & 
\\\noalign{\medskip} & &1& & & & & & & & \\
\noalign{\medskip}a_4&a_3&a_2&a_1&a_0&a_1&a_2&a_3&a_4& & \\
\noalign{\medskip} & &a_4&a_3&a_2&a_1&a_0&a_1&a_2&a_3&a_4\\
\noalign{\medskip}&&&\ddots&\ddots&\ddots&\ddots&\ddots&\ddots&\ddots&\ddots\\
\end {array} \right)	
\]

$a_{0-4}$ is listed in previous subsection.

\subsection{Sixth-Order Case}

Sixth order accurate case is given below. Notice that the stability condition for the upwind formulation could not be achieved. The discrete norm is given by $H=h\,diag(\frac{13649}{43200},\frac{12013}{8640},\frac{2711}{4320},\frac{5359}{4320},\frac{7877}{8640},\frac{43801}{43200},1,\ldots)$.

The upper left corner of $I_{F2C}$ is given by:
\[
{\scriptsize
\begin{array}{ll}
\begin{array}{l}
\medskip a_{1,1}={\frac {493203253466037198549455362364169948587300751}{1055237669662369450269932613530339432431112782}}\\
\medskip a_{1,2}= {\frac {2802735}{3494144}}\\
\medskip a_{1,3}= {\frac {65307081693884397225173945017702348750445730}{527618834831184725134966306765169716215556391}}\\
\medskip a_{1,4}= -{\frac {729595}{6988288}}\\
\medskip a_{1,5}= -{\frac {114866096841096884636391939697630554903763548}{527618834831184725134966306765169716215556391}}\\
\medskip a_{1,6}= -{\frac {1668147}{6988288}}\\
\medskip a_{1,7}= -{\frac {58150488726161516985950806005548507936924127}{527618834831184725134966306765169716215556391}}\\
\medskip a_{1,8}= {\frac {15525}{436768}}\\
\medskip a_{1,9}= {\frac {160596781068985961753206980048680817748946115}{1055237669662369450269932613530339432431112782}}\\
\medskip a_{1,10}= {\frac {89725}{436768}}\\
\medskip a_{1,11}= {\frac {28603010322006197876985672016642805153037153}{527618834831184725134966306765169716215556391}}\\
\medskip a_{1,12}= -{\frac {858195}{6988288}}\\
\medskip a_{1,13}= -{\frac {42754568564745224527531090194001857246708700}{527618834831184725134966306765169716215556391}}\\
\medskip a_{1,14}= {\frac {124485}{6988288}}\\
\medskip a_{1,15}= {\frac {10622334830569558938004175701898426132770338}{527618834831184725134966306765169716215556391}}\\
\medskip a_{1,16}= -{\frac {5409}{3494144}}\\
\medskip a_{1,17}= -{\frac {1157889404765086591410267452043339821491788}{527618834831184725134966306765169716215556391}}\\
 \medskip a_{2,1}= {\frac {10324674409544257975653283320299930288476692}{464377248357170642761107058624806491383799467}}\\
\medskip a_{2,2}= {\frac {1201479}{3075328}}\\
\medskip a_{2,3}= {\frac {63083470065854914766431215082154659242790665}{464377248357170642761107058624806491383799467}}\\
\medskip a_{2,4}= {\frac {721527}{3075328}}\\
\medskip a_{2,5}= {\frac {75656801964097450751057400512505712346209266}{464377248357170642761107058624806491383799467}}\\
\medskip a_{2,6}= {\frac {838645}{6150656}}\\
\medskip a_{2,7}= {\frac {30937093119827557151778923229206427591240999}{46437724835170642761107058624806491383799467}}\\
\medskip a_{2,8}= {\frac {345}{6150656}}\\
\end {array}
&
\begin{array}{l}
\medskip a_{2,9}=-{\frac {37294205185814962316045245318486095186507918}{
464377248357170642761107058624806491383799467}}\\
\medskip a_{2,10}=-{\frac {100125}{768832}}\\
\medskip a_{2,11}=-{\frac {36281493221658864449474225976906916398553827}{928754496714341285522214117249612982767598934}}\\
\medskip a_{2,12}={\frac {55219}{768832}}\\
\medskip a_{2,13}={\frac {22445907670907889073586224542480749734259694}{464377248357170642761107058624806491383799467}}\\
\medskip a_{2,14}=-{\frac {63387}{6150656}}\\
\medskip a_{2,15}=-{\frac {5433667362017414115320173550603825973179156}{464377248357170642761107058624806491383799467}}\\
\medskip a_{2,16}={\frac {5409}{6150656}}\\
\medskip a_{2,17}={\frac {578944702382543295705133726021669910745894}{464377248357170642761107058624806491383799467}}\\
\medskip a_{3,1}=-{\frac {5162337204772128987826641660149965144238346}{104797029908956098603626174638462532102012849}}\\
\medskip a_{3,2}=-{\frac {240439}{1388032}}\\
\medskip a_{3,3}={\frac {23386090748321137420688072323840400038565838}{104797029908956098603626174638462532102012849}}\\
\medskip a_{3,4}={\frac {105687}{173504}}\\
\medskip a_{3,5}={\frac {120607102321976279545796687898247496534070403}{419188119635824394414504698553850128408051396}}\\
\medskip a_{3,6}={\frac {75135}{694016}}\\
\medskip a_{3,7}=-{\frac {14476182846428389369733457155565337017487527}{209594059817912197207252349276925064204025698}}\\
\medskip a_{3,8}=-{\frac {217175}{1388032}}\\
\medskip a_{3,9}={\frac {6622670052963322216326122365241056731782031}{209594059817912197207252349276925064204025698}}\\
\medskip a_{3,10}={\frac {344085}{1388032}}\\
\medskip a_{3,11}={\frac {21321625226448423053859202901015977066304133}{209594059817912197207252349276925064204025698}}\\
\medskip a_{3,12}=-{\frac {78363}{694016}}\\
\medskip a_{3,13}=-{\frac {8305592718705970541713249560330845208737268}{104797029908956098603626174638462532102012849}}\\
\medskip a_{3,14}={\frac {5473}{347008}}\\
\medskip a_{3,15}={\frac {1892889085160894469318781649971017262255710}{104797029908956098603626174638462532102012849}}\\
\medskip a_{3,16}=-{\frac {1803}{1388032}}\\
\medskip a_{3,17}=-{\frac {192981567460847765235044575340556636915298}{104797029908956098603626174638462532102012849}}\\
\end {array}
\end {array}}
\]

The interior of $I_{F2C}$ is given by:
\[
\begin{array}{ll}
\begin{array}{l}
\medskip a_6={\frac {321100777805071156797079160300427016498}{38656226451108852306759931626138890483959}}\\
\medskip a_5={\frac {3}{512}}\\
\medskip a_4=-{\frac {1926604666830426940782474961802562098988}{38656226451108852306759931626138890483959}}\\
\medskip a_3=-{\frac {25}{512}}\\
\end {array}
&
\begin{array}{l}
\medskip a_2={\frac {4816511667076067351956187404506405247470}{38656226451108852306759931626138890483959}}\\
\medskip a_1={\frac {75}{256}}\\
\medskip a_0={\frac {25812195338906006034876765214121809824039}{77312452902217704613519863252277780967918}}\\
\end {array}
\end {array}
\]

\subsection{Sixth-Order Case, Non-SBP Interpolation}

The upper left part of $I_{C2F}$:
\[
\left( \begin {array}{cccccccccc} 1& & & & &&&&& \\
\noalign{\medskip}{\frac {63}{256}}&{\frac {315}{256}}&-{\frac {105}{128}}&\frac{63}{128}&-\frac{45}{256} &\frac{7}{256}&&&& \\
\noalign{\medskip} &1& & & & &&&&\\
\noalign{\medskip} -{\frac {7}{256}}&{\frac {105}{256}}&{\frac {105}{128}}&-\frac{35}{128}&\frac{21}{256} &-\frac{3}{256}&&&&\\
\noalign{\medskip} &&1 & & & &&&&\\
\noalign{\medskip}2a5& 2a_3&2a_1&2a_1&2a_3&2a_5&&&\\
\noalign{\medskip} 2a_6&2a_4&2a_2&2a_0&2a_2&2a_4&2a_6&&\\
\noalign{\medskip} & &2a5& 2a_3&2a_1&2a_1&2a_3&2a_5&\\
\noalign{\medskip} & &2a_6&2a_4&2a_2&2a_0&2a_2&2a_4&2a_6\\
\noalign{\medskip} &&&\ddots&\ddots&\ddots&\ddots&\ddots&\ddots&\\
\end{array}\right)
\]

The lower left part of $I_{F2C}$:
\[
\left( \begin {array}{ccccccccccccccc} 1& & & & & & & & & & &&&&
\\\noalign{\medskip} & &1& & & & & & & &&&&& \\
\\\noalign{\medskip} & && &1 & & & & & & &&&&\\
\noalign{\medskip}a_6&a_5&a_4&a_3&a_2&a_1&a_0&a_1&a_2&a_3&a_4&a_5 &a_6 \\
\noalign{\medskip} & &a_6&a_5&a_4&a_3&a_2&a_1&a_0&a_1&a_2&a_3&a_4&a_5 &a_6\\
\noalign{\medskip}&&&\ddots&\ddots&\ddots&\ddots&\ddots&\ddots&\ddots&\ddots&\ddots&\ddots&\ddots&\ddots\\
\end {array} \right)	
\]

$a_{0-6}$ is listed in previous subsection.

\subsection{Eighth-Order Case}

The upper left corner of $I_{F2C}$ is given by:

\[
{\scriptsize
\begin{array}{lll}
\begin{array}{l}
\medskip a_{1,1}=0.49670323986410765203\\
\medskip a_{1,2}= 0.80670591743192512511\\
\medskip a_{1,3}= -0.014476656476698073533\\
\medskip a_{1,4}= - 0.19497555250083395132\\
\medskip a_{1,5}= 0.16074268813765884450\\
\medskip a_{1,6}=  0.22100911221277072755\\
\medskip a_{1,7}=  -0.53042418939472250452\\
\medskip a_{1,8}=  - 0.86254139670456354517\\
\medskip a_{1,9}= 0.0064231825172866668299\\
\medskip a_{1,10}=  0.69727164011248914487\\
\medskip a_{1,11}= 0.61825742343094006838\\
\medskip a_{1,12}=  0.43307239695721032895\\
\medskip a_{1,13}=  -0.15848187861199173328\\
\medskip a_{1,14}=  - 0.83836831742920925562\\
\medskip a_{1,15}=  -0.41767239062238727391\\
\medskip a_{1,16}=   0.39252899650233765024\\
\medskip a_{1,17}= 0.33094736964907844809\\
\medskip a_{1,18}= - 0.073592865087013788440\\
 \medskip a_{1,19}=-0.099669762780958210515\\
\medskip a_{1,20}=  0.010748160731101219580\\
\medskip a_{1,21}=  0.018056833808814782786\\
 \medskip a_{1,22}= - 0.00077275234787125894193\\
 \medskip a_{1,23}= -0.0014911993994710636483\\
\medskip a_{2,1}= 0.0025487859584304862664\\
\medskip a_{2,2}=  0.47007080279424540800\\
\medskip a_{2,3}=   0.093589176759289320118\\
\medskip a_{2,4}=  0.33386224041716468423\\
\medskip a_{2,5}=  - 0.11418395501435496457\\
\medskip a_{2,6}=  -0.18998278818813064037\\
\medskip a_{2,7}=   0.38484431284491986603\\
\end {array}
&
\begin{array}{l}  
\medskip a_{2,8}=  0.65828018436035862232\\
\medskip a_{2,9}=   0.039704539050575728856\\
\medskip a_{2,10}=   -0.47252462545568042557\\
 \medskip a_{2,11}=  - 0.44892698918501097924\\
 \medskip a_{2,12}=  -0.34402935194949605213\\
\medskip a_{2,13}=    0.087284452472801643396\\
\medskip a_{2,14}=  0.59539401290875351190\\
\medskip a_{2,15}=   0.30484430526276345964\\
\medskip a_{2,16}=   -0.27096308216867871783\\
\medskip a_{2,17}=   - 0.23108785344796321535\\
\medskip a_{2,18}=  0.050608234918774219796\\
\medskip a_{2,19}=   0.068801842319351676214\\
\medskip a_{2,20}=   -0.0073322707316320135247\\
 \medskip a_{2,21}=  - 0.012333488857215226757\\
\medskip a_{2,22}=  0.00052275043402033040521\\
 \medskip a_{2,23}=  0.0010087644967132781708\\
\medskip a_{3,1}= - 0.022656973277393401539\\
\medskip a_{3,2}=- 0.93771184228725468159\\
\medskip a_{3,3}=- 0.052699680965389826267\\
\medskip a_{3,4}= 1.9581430555012001670\\
\medskip a_{3,5}= 1.6586148101136731602\\
 \medskip a_{3,6}= 2.7435837109255836264\\
\medskip a_{3,7}=  - 3.0164708462284474624\\
\medskip a_{3,8}=- 5.8235016129647971089\\
\medskip a_{3,9}=- 1.0472767830023645070\\
\medskip a_{3,10}= 3.2615209891555982371\\
\medskip a_{3,11}= 3.5478595826728208891\\
\medskip a_{3,12}=  3.0921640208240511054\\
\medskip a_{3,13}=  -0.34364793368678654581\\
\medskip a_{3,14}=  - 4.5566547247355317663\\
\medskip a_{3,15}=  -2.4329078834671892590\\
\end {array}
&
\begin{array}{l}  
\medskip a_{3,16}=   1.9882413967858906122\\
\medskip a_{3,17}=    1.7259601262271516763\\
\medskip a_{3,18}=    -0.36895458453625989435\\
\medskip a_{3,19}=    - 0.50448363278283993228\\
\medskip a_{3,20}=    0.052797763052066775846\\
\medskip a_{3,21}=     0.088972343374038343284\\
\medskip a_{3,22}=     -0.0037175165926095403240\\
\medskip a_{3,23}=     - 0.0071737841052106668688\\
\medskip a_{4,1}= 0.0021626748627943672418\\
\medskip a_{4,2}=   0.027636709246281031633\\
\medskip a_{4,3}= - 0.00055259484658035070394\\
\medskip a_{4,4}= - 0.033553649469391355855\\
\medskip a_{4,5}= - 0.049244245327794891233\\
\medskip a_{4,6}=  0.093489376222103426899\\
\medskip a_{4,7}=  0.45734620580442859300\\
 \medskip a_{4,8}=  0.66579609152388478842\\
 \medskip a_{4,9}= 0.24716623315221892947\\
 \medskip a_{4,10}= - 0.15893464065423852815\\
 \medskip a_{4,11}= -0.25881084331023040874\\
 \medskip a_{4,12}= - 0.26865808253702445366\\
\medskip a_{4,13}=  -0.018060020510041282529\\
\medskip a_{4,14}=   0.30841043055726240020\\
\medskip a_{4,15}=  0.17712461121229459766\\
\medskip a_{4,16}=  - 0.12549015561536110372\\
\medskip a_{4,17}=  -0.11250298507032009025\\
\medskip a_{4,18}=   0.022964117700293735857\\
\medskip a_{4,19}=  0.031709446610808019222\\
\medskip a_{4,20}=  - 0.0032149051440245356510\\
\medskip a_{4,21}=  -0.0054335286230388551025\\
\medskip a_{4,22}=   0.00022177994574852164638\\
\medskip a_{4,23}=  0.00042797426992744435493\\
\end {array}
\end {array}}
\]

 The interior of $I_{F2C}$ is given by:
\[
\begin{array}{ll}
\begin{array}{l}\medskip a_8=- 0.0023556211403911721282\\
\medskip a_7=-0.001220703125\\
\medskip a_6=0.018844969123129377026\\
\medskip a_5=0.011962890625\\
\medskip a_4=- 0.065957391930952819589\\
\end {array}
&
\begin{array}{l}  
\medskip a_3=- 0.059814453125\\
\medskip a_2=0.13191478386190563918\\
\medskip a_1=0.299072265625\\
\medskip a_0=0.33510652017261795103\\
\end {array}
\end {array}
\]

\subsection{Eighth- to Fourth-Order Case}

The upper left corner of $I_{F2C}$ is given by:

\[
{\scriptsize
\begin{array}{lll}
\begin{array}{l}
\medskip a_{1,1}=0.58483823651420605520\\
\medskip a_{1,2}= 0.65474734710702549940\\
\medskip a_{1,3}=- 0.083569673821672843574\\
 \medskip a_{1,4}=0.022425968859768418924\\
\medskip a_{1,5}= 0.31477269554917488790\\
\medskip a_{1,6}= - 0.19106598818802213021\\
\medskip a_{1,7}= -0.69007595018747002513\\
\medskip a_{1,8}= - 0.30500476635452414633\\
\medskip a_{1,9}= 0.011380303429615219398\\
\medskip a_{1,10}=  0.39023505096198663290\\
\medskip a_{1,11}= 0.69293443221288252711\\
\medskip a_{1,12}=  0.27930147563078413713\\
\medskip a_{1,13}=  -0.17509858064671991258\\
\medskip a_{1,14}=  - 0.36446581507281642455\\
\medskip a_{1,15}=  -0.53431810490188879575\\
\medskip a_{1,16}=  - 0.047645255801647499040\\
\medskip a_{1,17}=  0.48950384943343111452\\
\medskip a_{1,18}=   0.14968600350317559693\\
\medskip a_{1,19}=   -0.19713593239026365497\\
\medskip a_{1,20}=   - 0.051374514632076414476\\
\medskip a_{1,21}=   0.048757624195804968082\\
\medskip a_{1,22}=    0.0066939070511429668639\\
\medskip a_{1,23}=    -0.0055223124518961777586\\
\medskip a_{2,1}=    0.011443922859934377398\\
\medskip a_{2,2}= 0.37277008707234857179\\
\medskip a_{2,3}=  0.35282527696092662640\\
\medskip a_{2,4}= 0.096066986437174930272\\
\medskip a_{2,5}= - 0.22996810860646674853\\
\medskip a_{2,6}= 0.12123432064337482341\\
\medskip a_{2,7}=  0.49992261650980671917\\
\end {array}
&
\begin{array}{l}      
\medskip a_{2,8}= 0.25040596791141593267\\
\medskip a_{2,9}=  0.039631571628602904441\\
\medskip a_{2,10}=  -0.25355792154009083332\\
\medskip a_{2,11}=  - 0.50801440584689113247\\
\medskip a_{2,12}=  -0.21961026714024651967\\
\medskip a_{2,13}=   0.095341708111114624790\\
\medskip a_{2,14}=   0.24893964119173768791\\
\medskip a_{2,15}=    0.39235409374512417332\\
\medskip a_{2,16}= 0.042074861655253219822\\
\medskip a_{2,17}= - 0.34486664831340391286\\
\medskip a_{2,18}= -0.10508512775619684648\\
\medskip a_{2,19}=  0.13702066788245211707\\
\medskip a_{2,20}= 0.035285372237949024781\\
\medskip a_{2,21}= - 0.033422056958152836350\\
\medskip a_{2,22}= -0.0045282849361975826191\\
\medskip a_{2,23}=  0.0037357262504306790489\\
\medskip a_{3,1}=- 0.095787414722590617699\\
\medskip a_{3,2}=- 0.68323289165328697626\\
\medskip a_{3,3}=- 0.43569930844225202810\\
\medskip a_{3,4}= 1.8027331550089460412\\
\medskip a_{3,5}= 3.7567210660509388125\\
\medskip a_{3,6}= - 0.040511780395060609101\\
\medskip a_{3,7}= -3.9354784185565657354\\
\medskip a_{3,8}= - 2.5330463966461275528\\
\medskip a_{3,9}= - 1.0862631313184671074\\
\medskip a_{3,10}= 1.5476108406776833914\\
\medskip a_{3,11}=  4.0781324522987606857\\
\medskip a_{3,12}=   1.9733308525923243475\\
\medskip a_{3,13}=- 0.37565017560062195914\\
\medskip a_{3,14}=- 1.7659882182232900992\\
\medskip a_{3,15}=- 3.1675028227076565595\\
\end {array}
&
\begin{array}{l}  
\medskip a_{3,16}=- 0.41961367860804052329\\
\medskip a_{3,17}= 2.6138186452177357535\\
\medskip a_{3,18}= 0.79040581189389266572\\
\medskip a_{3,19}=- 1.0153714219277431776\\
\medskip a_{3,20}=- 0.25660116699213786317\\
\medskip a_{3,21}= 0.24235775601107244412\\
\medskip a_{3,22}= 0.032202698057872206652\\
\medskip a_{3,23}= - 0.026566452015385539718\\
\medskip a_{4,1}= 0.0085760788432811750313\\
\medskip a_{4,2}=  0.013512845571180982502\\
\medskip a_{4,3}= 0.022567253340913112846\\
\medskip a_{4,4}= - 0.058307759363885212745\\
\medskip a_{4,5}= -0.12041195448680664419\\
\medskip a_{4,6}=  0.22948717483369626712\\
\medskip a_{4,7}= 0.54863763561981341312\\
\medskip a_{4,8}=  0.41924435198770363142\\
\medskip a_{4,9}= 0.25222427818457038144\\
\medskip a_{4,10}= - 0.029868530356228792693\\
\medskip a_{4,11}= -0.30540634608037589943\\
\medskip a_{4,12}= - 0.18020289186124353787\\
\medskip a_{4,13}= -0.015353428793734303657\\
\medskip a_{4,14}=  0.10180220964078586702\\
\medskip a_{4,15}= 0.23447509646113431046\\
\medskip a_{4,16}=  0.040929410205375688302\\
\medskip a_{4,17}=  -0.17564595949066362349\\
\medskip a_{4,18}=  - 0.051894216445801530355\\
\medskip a_{4,19}=  0.065022578066260344054\\
\medskip a_{4,20}=   0.015872184741392310993\\
\medskip a_{4,21}=   -0.014923762934166461295\\
\medskip a_{4,22}=   - 0.0019211515134671124718\\
\medskip a_{4,23}= 0.0015849038302656338915\\
\end {array}
\end {array}}
\]

 The interior of $I_{F2C}$ is given by:
\[
\begin{array}{ll}
\begin{array}{l}
\medskip a_8= - 0.0087234986549392520415\\
\medskip a_7= 0.010574245783011317441\\
\medskip a_6= 0.041177278191813548941\\
\medskip a_5=- 0.037706537233780062173\\
\medskip a_4=- 0.090162662747982636209\\
\end {array}
&
\begin{array}{l} 
\medskip a_3= 0.018424137003272281874\\
\medskip a_2= 0.12963112574463663703\\
\medskip a_1= 0.25870815444749646286\\
\medskip a_0= 0.35615551493294340455\\
\end {array}
\end {array}
\]

\bibliographystyle{plain}
\bibliography{referenser}

\begin{thebibliography}{10}

\bibitem{Abarbanel97}
S.~Abarbanel and A.~Ditkowski.
\newblock Asymptotically stable fourth-order accurate schemes for the diffusion
  equation on complex shapes.
\newblock {\em J. Comput. Phys.}, 133:279--288, 1997.

\bibitem{Bayliss86}
A.~Bayliss, K.~E. Jordan, B.~J. Lemesurier, and E.~Turkel.
\newblock A fourth order accurate finite difference scheme for the computation
  of elastic waves.
\newblock {\em Bull. Seismol. Soc. Amer.}, 76(4):1115--1132, 1986.

\bibitem{CGA}
M.~H. Carpenter, D.~Gottlieb, and S.~Abarbanel.
\newblock {T}ime-stable boundary conditions for finite-difference schemes
  solving hyperbolic systems: {M}ethodology and application to high-order
  compact schemes.
\newblock {\em J. Comput. Phys.}, 111(2), 1994.

\bibitem{CarpenterNordstrom99}
M.~H. Carpenter, J.~Nordstr\"om, and D.~Gottlieb.
\newblock A {S}table and {C}onservative {I}nterface {T}reatment of {A}rbitrary
  {S}patial {A}ccuracy.
\newblock {\em J. Comput. Phys.}, 148, 1999.

\bibitem{Diener_07}
Peter Diener, Ernst~Nils Dorband, Erik Schnetter, and Manuel Tiglio.
\newblock Optimized high-order derivative and dissipation operators satisfying
  summation by parts, and applications in three-dimensional multi-block
  evolutions.
\newblock {\em J. Sci. Comput.}, 32(1):109--145, 2007.

\bibitem{Gustafsson81}
B.~Gustafsson.
\newblock The convergence rate for difference approximations to general mixed
  initial boundary value problems.
\newblock {\em SIAM J. Numer. Anal.}, 18(2):179--190, Apr. 1981.

\bibitem{GustafssonKreissOliger}
B.~Gustafsson, H.-O. Kreiss, and J.~Oliger.
\newblock {\em Time dependent problems and difference methods}.
\newblock John Wiley \& Sons, Inc., 1995.

\bibitem{GKS}
B.~Gustafsson, H.~O. Kreiss, and A.~Sundstr\"om.
\newblock Stability theory of difference approximations for mixed initial
  boundary value problems.
\newblock {\em Math. Comp.}, 26(119), 1972.

\bibitem{bertil_pelle2}
B.~Gustafsson and P.~Olsson.
\newblock Fourth-order difference methods for hyperbolic {IBVPs}.
\newblock {\em J. Comput. Physics}, 117(1), 1995.

\bibitem{Hesthaven98}
J.S. Hesthaven.
\newblock A stable penalty method for the compressible {N}avier-{S}tokes
  equations: {II} {M}ulti-dimensional domain decomposition schmemes.
\newblock {\em SIAM {J}.{S}ci.{C}omput.}, 1998.

\bibitem{kopriva96}
D.~A. Kopriva and J.~H. Kolias.
\newblock A conservative staggered-grid chebyshev multidomain method for
  compressible flows.
\newblock {\em J. Comput. Phys.}, 125(1):244--261, 1996.

\bibitem{KreissScherer74}
H.-O. Kreiss and G.~Scherer.
\newblock Finite element and finite difference methods for hyperbolic partial
  differential equations.
\newblock {\em Mathematical Aspects of Finite Elements in Partial Differential
  Equations., Academic Press, Inc.}, 1974.

\bibitem{KreissOliger72}
Heinz-Otto Kreiss and Joseph Oliger.
\newblock Comparison of accurate methods for the integration of hyperbolic
  equations.
\newblock {\em Tellus XXIV}, 3, 1972.

\bibitem{Lehner_05}
L.~Lehner, O.~Reula, and M.~Tiglio.
\newblock Multi-block simulations in general relativity: high-order
  discretizations, numerical stability and applications.
\newblock {\em Classical Quantum Gravity}, 22:5283--5321, 2005.

\bibitem{kn:Lele(1992)}
S.~K. Lele.
\newblock Compact finite difference schemes with spectral-like resolution.
\newblock {\em J. Comput. Phys.}, 103:16--42, 1992.

\bibitem{Mattsson03}
K.~Mattsson.
\newblock Boundary procedures for summation-by-parts operators.
\newblock {\em Journal of Scientific Computing}, 18:133--153, 2003.

\bibitem{MattssonHam08}
K.~Mattsson, F.~Ham, and G.~Iaccarino.
\newblock Stable and accurate wave propagation in discontinuous media.
\newblock {\em J. Comput. Phys.}, In Press, 2008.

\bibitem{MattssonNordstrom04}
K.~Mattsson and J.~Nordstr\"om.
\newblock Summation by parts operators for finite difference approximations of
  second derivatives.
\newblock {\em J. Comput. Phys.}, 199(2):503--540, 2004.

\bibitem{MattssonNordstrom06}
K.~Mattsson and J.~Nordstr\"om.
\newblock High order finite difference methods for wave propagation in
  discontinuous media.
\newblock {\em J. Comput. Phys.}, 220:249--269, 2006.

\bibitem{MattssonSvard06}
K.~Mattsson, M.~Sv\"ard, M.H. Carpenter, and J.~Nordstr\"om.
\newblock High-order accurate computations for unsteady aerodynamics.
\newblock {\em Computers \& Fluids}, 36:636--649, 2006.

\bibitem{MattssonSvard04}
K.~Mattsson, M.~Sv\"ard, and J.~Nordstr\"om.
\newblock {S}table and {A}ccurate {A}rtificial {D}issipation.
\newblock {\em Journal of Scientific Computing}, 21(1):57--79, August 2004.

\bibitem{MattssonSvard07}
K.~Mattsson, M.~Sv\"ard, and M.~Shoeybi.
\newblock Stable and accurate schemes for the compressible navier-stokes
  equations.
\newblock {\em J. Comput. Phys.}, 227(4):2293--2316, 2008.

\bibitem{NordstromCarpenter99}
J.~Nordstr\"om and M.~H. Carpenter.
\newblock Boundary and interface conditions for high-order finite-difference
  methods applied to the {E}uler and {N}avier-{S}tokes equations.
\newblock {\em J. Comput. Phys.}, 148:341--365, 1999.

\bibitem{NordstromCarpenter01}
J.~Nordstr\"om and M.~H. Carpenter.
\newblock High-order finite difference methods, multidimensional linear
  problems, and curvilinear coordinates.
\newblock {\em J. Comput. Phys.}, 173:149--174, 2001.

\bibitem{Nordstrom_gong_2006}
J.~Nordstr\"om and J.~Gong.
\newblock A stable hybrid method for hyperbolic problems.
\newblock {\em J. Comput. Phys.}, 212:436--453, 2006.

\bibitem{NordstromMattsson07}
J.~Nordstr\"om, K.~Mattsson, and R.C. Swanson.
\newblock Boundary conditions for a divergence free velocity-pressure
  formulation of the incompressible navier-stokes equations.
\newblock {\em J. Comput. Phys.}, 225:874--890, 2007.

\bibitem{Olsson95(1)}
P.~Olsson.
\newblock Summation by parts, projections, and stability {I}.
\newblock {\em Math. Comp.}, 64:1035, 1995.

\bibitem{Olsson95(2)}
P.~Olsson.
\newblock Summation by parts, projections, and stability {II}.
\newblock {\em Math. Comp.}, 64:1473, 1995.

\bibitem{DeRangoZingg01}
S.~De Rango and D.~W. Zingg.
\newblock High-order aerodynamic computations on multi block grids.
\newblock {\em AIAA Paper}, 2001-2631, 2001.

\bibitem{Sjogreen_95}
B.~Sjogreen.
\newblock High order centered difference methods for the compressible
  navier-stokes equations.
\newblock Technical report 01.01, RIACS, NASA Ames Research Center, 2001.

\bibitem{Strand94}
B.~Strand.
\newblock Summation by parts for finite difference approximations for d/dx.
\newblock {\em J. Comput. Physics}, 110:47--67, 1994.

\bibitem{Strikwerda97}
John~C. Strikwerda.
\newblock High-order-accurate schemes for incompressible viscous flow.
\newblock {\em International Journal for Numerical Methods in Fluids},
  24:715--734, 1997.

\bibitem{SvardNordstrom08a}
M.~Sv\"ard, M.~H. Carpenter, and J.~Nordstr\"om.
\newblock A stable high-order finite difference scheme for the compressible
  {N}avier--{S}tokes equations, far-field boundary conditions.
\newblock {\em J. Comput. Physics}, 225:1020--1038, February 2008.

\bibitem{SvardNordstrom08b}
M.~Sv\"ard, M.~H. Carpenter, and J.~Nordstr\"om.
\newblock A stable high-order finite difference scheme for the compressible
  {N}avier--{S}tokes equations, no-slip wall boundary conditions.
\newblock {\em J. Comput. Physics}, 227:4805--4824, May 2008.

\bibitem{Svard04_Mattsson}
M.~Sv\"ard, K.~Mattsson, and J.~Nordstr\"om.
\newblock Steady-state computations using summation-by-parts operators.
\newblock {\em Journal of Scientific Computing}, 24(1):79--95, July 2005.

\bibitem{SvardNordstrom06}
M.~Sv\"ard and J.~Nordstr\"om.
\newblock On the order of accuracy for difference approximations of
  initial-boundary value problems.
\newblock {\em J. Comput. Physics}, 218:333--352, October 2006.

\bibitem{ZinggDeRango99}
D.~W. Zingg, S.~De Rango, M.~Nemec, and T.~H. Pulliam.
\newblock Comparison of several spatial discretizations for the
  {N}avier-{S}tokes equations.
\newblock {\em AIAA Paper}, 99-3260, 1999.

\end{thebibliography}

\end{document}